\documentclass[a4paper, 11pt]{amsart}
\usepackage[latin1]{inputenc}
\usepackage[ T1]{fontenc}
\usepackage[english]{babel}
\usepackage{amssymb}
\usepackage{amsmath}
\usepackage{amsthm}
\usepackage{amscd}
\usepackage{amsfonts}
\usepackage{stmaryrd}
\usepackage{pb-diagram}
\usepackage{epic,eepic,epsfig}
\usepackage{a4wide}
\usepackage{nextpage}
\usepackage{fancyhdr}
\usepackage{enumerate}
\usepackage{hyperref}
\usepackage{float}
\usepackage{tikz}

\newtheorem{prop}{Proposition}[section]
\newtheorem{corollary}[prop]{Corollary}

\newtheorem{lemma}[prop]{Lemma}
\newtheorem{thm}[prop]{Theorem}

\theoremstyle{definition}
\newtheorem{rem}[prop]{Remark}

\newtheorem{definition}[prop]{Definition}

\newcommand{\GL} {\mathop{\mathrm{GL}}}
\newcommand{\PGL} {\mathop{\mathrm{PGL}}}

\newcommand{\ti} {\text{i}}
\newcommand{\Tcal}{\mathcal{T}}
\newcommand{\Ocal}{\mathcal{O}}
\newcommand{\Ifra}{\mathfrak{I}}
\newcommand{\hTag}{h_{\mathrm{Tag}}}
\newcommand{\Fm}{\mathfrak{m}}

\def\BC{\mathbb{C}}
\def\BF{\mathbb{F}}
\def\BQ{\mathbb{Q}}
\def\BZ{\mathbb{Z}}

\def\BR{\mathbb{R}}

\def\CF{\mathcal{F}}
\newcommand{\CO}{\mathcal{O}}

\def\Cinf{\mathbb{C}_\infty}

\def\GL{\mathrm{GL}}

\def\Tag{\mathrm{Tag}}

\usepackage[OT2,T1]{fontenc}
\DeclareSymbolFont{cyrletters}{OT2}{wncyr}{m}{n}
\DeclareMathSymbol{\Sha}{\mathalpha}{cyrletters}{"58}

\begin{document}

\title{Heights of Drinfeld modular polynomials and Hecke images}
\author{}
\author{Florian Breuer \and Fabien Pazuki \and Zhenlin Ran}

\address{School of Computer and Information Sciences, The University of Newcastle,
University Drive, Callaghan, NSW 2308, Australia.}
 \email{Florian.Breuer@newcastle.edu.au}

\address{Department of Mathematical Sciences, University of Copenhagen,
 Universitetsparken 5, 
2100 Copenhagen \O, Denmark, and Universit\'e de Bordeaux, 33405 Talence, France.}
 \email{fpazuki@math.ku.dk}

 \address{Department of Mathematics, The University of Hong Kong, Run Run Shaw Building, Pokfulam, Hong Kong}
 \email{zhenlinran@gmail.com}

\maketitle

\noindent \textbf{Abstract.}
We obtain explicit upper and lower bounds on the size of the coefficients of the Drinfeld modular polynomials $\Phi_N$ for any monic $N\in\BF_q[t]$.
These polynomials vanish at pairs of $j$-invariants of Drinfeld $\BF_q[t]$-modules of rank 2 linked by cyclic isogenies of degree $N$.
The main term in both bounds is asymptotically optimal as $\deg N$ tends to infinity. We also obtain precise estimates on the Weil height and Taguchi height of Hecke images of Drinfeld modules of rank~2.

{\flushleft
\textbf{Data Availability Statement:} There is no data associated to this text.\\
\textbf{Conflict of Interest Statement:} The authors declare no conflicts of interest.
}

%11F52: Modular forms associated to Drinfeld modules
%11G09: Drinfeld modules, higher dimensional motives
%11R58: Arithmetic theory of algebraic function fields
%14H25: Arithmetic ground fields for curves

\begin{center}
---------
\end{center}

\thispagestyle{empty}

\maketitle

\section{Introduction}\label{section def}

Let $A = \BF_q[t]$ and $F = \BF_q(t)$, equipped with the absolute value $\vert\cdot\vert$ given by $|x|=q^{\deg x}$ for any $x\in{F}$, corresponding to the place $\infty$, where $\deg$ is the degree in $t$ for a rational fraction in $t$. We denote by $F_\infty = \BF_q((\frac{1}{t}))$ the completion of $F$ at $\infty$ and by $\Cinf$ the completion of an algebraic closure of $F_\infty$. The absolute value $|\cdot|$ can be naturally extended to $\Cinf$, which we again denote by $|\cdot|$. Let us introduce the non-commutative ring of twisted polynomials \[F\{\tau\}=\left\{ \sum_{i=0}^n c_i \tau^i \,\Big|\, n \geq 0, c_i \in F\right\},\]
with the classical additive law, and the multiplication defined by $\tau c= c^q \tau$ for every $c \in F$.

Let $\phi$ be a Drinfeld $A$-module of rank $2$ over $\Cinf$. It is given by a twisted polynomial in $\tau$, of the form
$$ \phi_t = t + g \tau + \Delta \tau^2,$$
where $g$ and $\Delta$ are in $\BC_{\infty}$ and $\Delta \neq 0$. The $j$-invariant of $\phi$ is an element of $\BC_{\infty}$ defined as
$$ j(\phi) = \frac{g^{q+1}}{\Delta},$$
and it characterises the isomorphism class of the Drinfeld module $\phi$ over $\Cinf$. We will only be concerned with Drinfeld modules of rank 2 in this work.

Let $N\in A$ be a monic polynomial. 
Analoguously to the case of elliptic modular polynomials, S.~Bae \cite{Bae92} defined the {\em Drinfeld modular polynomial} $\Phi_N(X,Y)\in A[X,Y]$ (see also Definition \ref{drmopoly} below).

It vanishes precisely at the pairs of $j$-invariants $(j(\phi), j(\phi'))$, where $\phi$ and $\phi'$ are rank 2 Drinfeld modules linked by an isogeny $f : \phi \rightarrow \phi'$ with kernel $\ker f \cong A/NA$. For $\deg N > 0$, this polynomial is symmetric in $X$ and $Y$ and its degree in each variable equals
\[
\psi(N) = |N|\prod_{P|N}\left(1 + \frac{1}{|P|}\right).
\]

This polynomial defines an integral model of the Drinfeld modular curve $X_0(N)$. Such modular curves are interesting for several reasons. They are used for constructing towers of curves with many rational points, with applications in coding theory, see for example \cite{Elkies}. These curves also parametrise elliptic curves over function fields via a theorem in analogy with the Taniyama-Shimura-Weil (former) conjecture, see \cite{GeRe}. More precisely let $E$ be a non-isotrivial elliptic curve over $F$, with split multiplicative reduction at the place $\infty$. Then there exists a modular parametrisation $X_0(N)\rightarrow E$. See \cite{Petit} for recent explicit computations. 

The coefficients of Drinfeld modular polynomials grow rapidly with $|N|$. This is in clear analogy with the case of  classical modular polynomials, which are built using cyclic isogenies of fixed degree between elliptic curves.

Our main goal is to estimate the size of the coefficients of $\Phi_N$, \textit{i.e.} the height of $\Phi_N$. We define the {\em height} of a non-zero polynomial $f$ with coefficients in $\Cinf$ to be 
\[
h(f) := \log \max_c |c|,
\]
where $c$ ranges over all coefficients of $f$, and $\log=\log_q$ is the logarithm in base $q$.

There exist some results in this direction. Hsia proved in \cite{hsia} that, as $\deg N \rightarrow \infty$, the height of $\Phi_N$ satisfies
\[
h(\Phi_N) = \frac{q^2-1}{2}\psi(N)\left(\deg N - 2\kappa_N + O(1)\right),
\]
where
\[
\kappa_N := \sum_{P|N}\frac{\deg P}{|P|},
\]
but the bounded term $O(1)$ was not made explicit. Explicit upper bounds on $h(\Phi_N)$ appear in \cite{BL97} and \cite{BPR}, but neither of these results are asymptotically optimal.

In the present paper, we obtain an explicit, asymptotically optimal estimate for $h(\Phi_N)$ and derive some consequences of our methods. Our first result is the following theorem. Define

 \begin{equation}
    \lambda_N := \sum_{P^n \| N}\frac{|P|^n-1}{|P|^{n-1}(|P|^2-1)}\deg P,
 \end{equation}
then we have the following result.

\begin{thm}\label{thm:main}
    Let $N\in A$ be monic. Then the height of the modular polynomial $\Phi_N(X,Y)$ satisfies
    \[
    h(\Phi_N) = \frac{q^2-1}{2}\psi(N)\left(\deg N - 2\lambda_N + b_q(N)\right),
    \]
    where
    \[
    0\leq b_q(N) \leq 4 + \frac{5q^3-q^2-q+1}{q(q+1)(q-1)^2}.
    \]
\end{thm}

The proof is given in Section \ref{proof of main}. The approach is similar, at least partially, to that in \cite{BGP} and \cite{BP}, where analogous results are proved for elliptic modular polynomials. The general idea is to estimate the Mahler measure of $\Phi_N(X,j(z))$ for some chosen $z$. Then using the result that the geometric mean of the Hecke image of the discriminant function equals  the discriminant function itself, we can reduce our proof to estimating the sum of the logarithms of ``imaginary parts'' (see equation (\ref{imaginary part}) for the definition) of the representatives of the Hecke images of the corresponding Drinfeld modules on the upper half plane. Function fields Farey sequences are used for these estimates. This function field analogue of Farey sequences is worked out in \cite[Section 3]{hsia}, see also below in paragraph \ref{Farey}.

However, there are new features in function field arithmetic. One new feature of the current paper is that we show $h(\Phi_N(X,Y)) = h(\Phi_N(X,y))$ for a single $y\in\overline{\BF_q}$, whereas in previous papers (for instance in \cite{Paz}) one needed to specialize at a number of $Y$-values and then use Lagrange interpolation to estimate $h(\Phi_N)$. This trick relies on the non-Archimedean metric on $\Cinf$. Another new feature concerns the upper half plane. Unlike the classical case where a good geometry exists in the Poincar\'e upper half plane, in this paper it is not straightforward to control some ``large'' representatives from the Hecke images. The Bruhat-Tits tree comes into play to address this issue. Concretely, we control the fundamental representatives of given vertices on the Bruhat-Tits tree (cf. Theorem \ref{fdmr}), which helps us in our estimates. 

The upper bound in Theorem \ref{thm:main} is probably not sharp, in particular we expect the constant $4$ to be too large, but the method in place doesn't yield better.
\\

As another application of our techniques, we compare the height of a Drinfeld module of rank 2 with the heights of its images under a Hecke correspondence. The results may be viewed as ``Hecke-averaged'' versions of the results in \cite{BPR}.

There are at least two notions of height of a Drinfeld module of rank 2 that one may consider: the Weil height of the $j$-invariant, and the stable Taguchi height, see Section \ref{heights} for definitions. 

In our second main result, we prove the following analog
of a result of Silverman \cite{Sil}, comparing the Weil height of a $j$-invariant with the height of its Hecke image.

\begin{thm}\label{Silverman}
    Let $\phi$ be a Drinfeld module of rank 2 defined over a global function field $K$. Then one has
    \begin{align*}
        % -\frac{1}{\psi(N)}h(\Phi_N)  
        \frac{q^2-1}{2}\big(2\lambda_N - \deg N - b_q(N)\big)
        & \leq \;
        h(j(\phi)) - \frac{1}{\psi(N)}\sum_{C}h(j(\phi/C)) \\
        & \leq q +\frac{q^2-1}{2}\min \left[0, \;
        2\lambda_N - \deg N + \log\left(\frac{1}{q}h(j(\phi)) + 1\right)
        \right],
    \end{align*}
    where $C$ ranges over all submodules $C\subset \phi[N]$ isomorphic to $A/NA$ and
    $b_q(N)$ is the bounded term from Theorem \ref{thm:main}.
\end{thm}

As a third result, we prove the following analog of a result of Autissier \cite{Aut}, which shows that the stable Taguchi height gives us much cleaner expressions than the Weil height of the $j$-invariant when estimating the heights of Hecke images. We denote by $h_\Tag(\phi)$ the stable Taguchi height of $\phi$ (see Section \ref{heights}, or \cite{Tag}, for definitions).

\begin{thm}\label{Autissier}
    Let $\phi$ be a Drinfeld module of rank 2 defined over a global function field $K$, and suppose that $\phi$ and $\phi/C$ have everywhere stable reduction over $K$, for all submodules $C\subset \phi[N]$ with $C\cong A/NA$. Then
    \[
    h_\Tag(\phi) - \frac{1}{\psi(N)} \sum_{C}h_\Tag(\phi/C) = 
    \lambda_N - \frac{1}{2}\deg N.
    \]
\end{thm}

The proofs are given in Section \ref{proof of Hecke}. We decompose the heights into a part coming from finite places (treated using Tate uniformisation) and a part coming from the analytic uniformisation of Drinfeld modules. We then make use of transformation properties of the Drinfeld modular discriminant proved in Proposition \ref{davg}. We conclude by an application of Lemma \ref{imagsum}, which is inspired by previous work of Autissier \cite{Aut} in characteristic zero.

\section*{Acknowledgements}
The authors thank the IRN GandA (CNRS) for support. FP is supported by ANR-20-CE40-0003 Jinvariant. FB is supported by the Alexander-von-Humboldt Foundation.

\section{Setting}

We gather in this section several crucial tools which will be used in the proofs of 
our main results.

We start with some background in height theory, then we cover results on the Drinfeld fundamental domain, Bruhat-Tits trees, the building map, Drinfeld modular forms, and Farey sequences in the polynomial setting. The last paragraph builds on the previous ones and provides the proof of an estimate on the size of some representatives in the fundamental domain. 

\subsection{Heights}\label{heights}
We consider two notions of height of a Drinfeld module. First, the Weil height of the $j$-invariant and second, the stable Taguchi height. 

Let $K/F$ be a finite extension and 
denote by $h(x)$ the usual Weil height for an element $x\in K$:
\[
h(x) := \frac{1}{[K:F]}\sum_{v\in M_K} n_v \log\max(1, |x|_v),
\]
where $M_K$ denotes the set of places of $K$, $n_v = [K_v : F_w]$ is the local degree with $v|w$ and $w\in M_F$, and $|x|_v$ is the absolute value associated to the place $v\in M_K$ normalised as follows:
\[
|x|_v:=|N_{K_v/F_w}(x)|_w^{1/n_v},
\]
where $N_{K_v/F_w}(\cdot)$ is the norm. The absolute values and local degrees have the following properties:
\[
\sum_{v\in M_K} n_v\log|x|_v = 0 \quad\text{(product formula for $x\neq0$), and}\]
\[
\sum_{v|w}n_v = [K:F]\quad\text{for any $w\in M_F$} \quad\text{(extension formula).}
\]
We also denote by $M_K^\infty$ the set of infinite places, \textit{i.e.} those $v|\infty$, and by $M_K^f = M_K \smallsetminus M_K^\infty$, the set of finite places. The infinite places $v\in M_K^\infty$ may be represented by $|x|_v = |\sigma(x)|$ for any $x\in{K}$, where $\sigma : K \hookrightarrow \Cinf$ is an embedding of $K$ associated to $v$. 

Let $\phi$ be a Drinfeld module defined over a finite extension  $K/F$. We recall that $\phi$ is said to have {\em stable reduction} at a place $v\in M_K^f$ if it is isomorphic to a Drinfeld module $\tilde\phi$ defined over the valuation ring $\CO_v$ of $v$ whose reduction modulo the maximal ideal $\Fm_v$ of $\CO_v$ is a Drinfeld module of positive rank over the residue field $\CO_v/\Fm_v$. We say that $\phi$ has {\em everywhere stable reduction} if it has stable reduction at every finite place $v\in M_K^f$.
By \cite[Lemme 2.10]{DD}, every Drinfeld module over $K$ acquires everywhere stable reduction after replacing $K$ by a finite extension thereof.

In \cite{Tag} Taguchi defines the {\em differential height} of $\phi$ as the Arakelov degree of the metrised conormal line-bundle along the unit section associated to a minimal model of~$\phi$. It serves as the analogue of the Faltings height of an abelian variety. 
Let us invoke the identity (5.9.1) of \cite{Tag}, valid for Drinfeld modules of rank $r$ with everywhere stable reduction, which we adopt as our definition (with $r=2$):
\begin{align} \label{eq:Tag height}
    \hTag(\phi) = &  \frac{1}{[K:F]} \sum_{v\in M_K^f} n_v\log \max\left(|g(\phi)|_v^{1/(q-1)},|\Delta(\phi)|_v^{1/(q^2-1)}\right)\\ \nonumber
    & - \frac{1}{[K:F]}\sum_{v\in M_K^\infty} n_v \log D(\Lambda_v)^{1/2}. 
\end{align}

% \begin{equation}
%    \hTag(\phi) = \hTag^f(\phi) + \hTag^\infty(\phi),
% \end{equation}

% \noindent where we pose
% \begin{align*}
%     \hTag^f(\phi) & :=\frac{1}{[K:F]} \sum_{v\in M_K^f} n_v\log \max(|g(\phi)|_v^{1/(q-1)},|\Delta(\phi)|_v^{1/(q^2-1)}),\\
%     \hTag^\infty(\phi) & :=- \frac{1}{[K:F]}\sum_{v\in M_K^\infty} n_v \log D(\Lambda_v)^{1/2}.
% \end{align*}
Here, $\Lambda_v \subset\Cinf$ is the lattice which uniformises the Drinfeld module $\phi^v$ over $\Cinf$, given by $\phi_t^v = t + \sigma(g(\phi))\tau + \sigma(\Delta(\phi))\tau^2$, where $\sigma : K \hookrightarrow \Cinf$ is an embedding associated to the infinite place $v$; and $D(\Lambda_v)$ is its covolume, which we define in the next paragraph, see Definition \ref{covolume}.

Using the definition $j(\phi) = g(\phi)^{q+1}/\Delta(\phi)$ and the product formula applied to $\Delta(\phi)$, we rewrite (\ref{eq:Tag height}) in the following form, which has the advantage that each local term is an isomorphism invariant:
\begin{align}\label{eq: Tag height 2}
    \hTag(\phi) = &  \frac{1}{(q^2-1)[K:F]} \sum_{v\in M_K^f} n_v\log \max\left(|j(\phi)|_v, 1\right) \\ \nonumber 
    & - \frac{1}{(q^2-1)[K:F]} \sum_{v\in M_K^\infty} n_v \log \left(|\Delta(\phi)|_v D(\Lambda_v)^{(q^2-1)/2}\right). 
\end{align}

% \[
% \hTag(\phi) =   \frac{1}{(q^2-1)[K:F]}\left[ \sum_{v\in M_K^f} n_v\log \max\left(|j(\phi)|_v, 1\right) - \sum_{v\in M_K^\infty} n_v \log \left(|\Delta(\phi)|_v D(\Lambda_v)^{(q^2-1)/2}\right)\right]
% \]

\subsection{Lattices and the Drinfeld upper half plane}

Recall that an $A$-\textit{lattice in} $\Cinf$ is a discrete and finitely generated $A$-submodule $\Lambda\subset \Cinf$. In the sequel, we shall always refer to these $A$-lattices as ``lattices'' (with ``$A$'' omitted) and all our lattices are in $\Cinf$. The \textit{rank} of a lattice $\Lambda$ is the dimension of $\Lambda\otimes_A F$ as an $F$-vector space. 

A \textit{successive minimum basis} for a lattice $\Lambda$ of rank $r$ is a tuple $(\omega_1,...,\omega_r)$ such that $\omega_i\in \Lambda$ for each $i=1,...,r$ and such that
$$|\omega_1|=\inf\{|\lambda|:\lambda\in \Lambda\smallsetminus\{0\}\},\ \text{and }|\omega_{i+1}|=\inf \left\{|\lambda|:\lambda\in \Lambda \smallsetminus \bigoplus_{j=1}^i A\cdot \omega_j \right\}.$$
For any lattice, there always exists at least one successive minimum basis.
Any successive minimum basis $(\omega_1,...,\omega_r)$ in $\Lambda$ is an orthogonal basis (cf. \cite[Lemma 4.2]{Tag}) in the sense that
$$\Lambda=\bigoplus_{i=1}^r A\cdot \omega_i$$
and for any $\alpha_i\in F_{\infty}$,
$$\left|\sum_{i=1}^r \alpha_i \omega_i\right|=\max\{|\alpha_i\omega_i|: 1\leq i\leq r\}.$$

\begin{definition}\label{covolume}
        Let $\Lambda\subset \Cinf$ be a lattice of rank $r$ with a successive minimum basis $(\omega_1,...,\omega_r)$. The \textit{covolume} $D(\Lambda)$ is defined to be
        \[
        D(\Lambda):=|\omega_1|\cdots|\omega_r|.
        \]
\end{definition}

It is known (cf. \cite[\S4]{Tag}) that the magnitudes $|\omega_1|, \ldots, |\omega_r|$, and thus also $D(\Lambda)$, do not depend on the choice of successive minimal basis of $\Lambda$.

\medskip

Let $\Omega = \Cinf \smallsetminus F_\infty$ denote the {\em Drinfeld upper half plane}. Then every $z\in\Omega$ gives rise to a rank 2 lattice $\Lambda_z = zA + A \subset\Cinf$.

The group $\GL_2(F_\infty)$ acts on $\Omega$ via fractional linear transformations as follows. For $\gamma=\begin{pmatrix}
    a & b \\ c & d
\end{pmatrix}\in\GL_2(F_\infty)$ and $z\in\Omega$, we have
\[
\gamma(z) = \frac{az+b}{cz+d}.
\]
For $z\in\Omega$ we define
\begin{equation}\label{imaginary part}
    |z|_\ti := \inf_{x\in F_\infty}|z-x|,
\end{equation}
which one can think of as the size of the ``imaginary part'' of $z$. For 
$\gamma=\begin{pmatrix}
    a & b \\ c & d
\end{pmatrix}\in\GL_2(F_\infty)$ we have (see \textit{e.g.} \cite[(3.2)]{GekTur})
\begin{equation}\label{imtrans}
|\gamma(z)|_\ti = \frac{|\det(\gamma)|}{|cz+d|^2}|z|_\ti.
\end{equation}

Throughout this paper, we set $\Gamma = \GL_2(A)$. 
The set
\[
\mathcal{F}:=\{z\in \Omega:|z|=|z|_{\ti}\geq 1\},
\]
is called the {\em fundamental domain} in $\Omega$, due to the following property, see \textit{e.g.} \cite[(6.7)]{GekTur}.

\begin{prop}\label{fundamenta domain}
    \begin{enumerate}[1.]
        \item Each $z\in\Omega$ is $\Gamma$-equivalent to an element in $\CF$.
        \item Each $z\in\Omega$ is $\Gamma$-equivalent to only finitely many elements in $\CF$.
        \item If $z, z' \in\CF$ are $\Gamma$-equivalent, then $|z| = |z|_\ti = |z'| = |z'|_\ti$.
    \end{enumerate}
\end{prop}

We say that a lattice $\Lambda\subset\Cinf$ is {\em reduced} if $1\in \Lambda$ and $|\lambda|\geq 1$ for all $\lambda\in \Lambda\smallsetminus \{0\}.$ 

\begin{lemma}\label{smb}
    Let $\Lambda\subset \Cinf$ be a lattice of rank $r$ with $r\geq 2$. Then $\Lambda$ is reduced if and only if there exists a successive minimum basis $(\omega_1,\omega_2,...,\omega_r)$ such that $\omega_1=1$ and $\omega_j\in \mathcal{F}$ for all $2\leq j\leq r$, and moreover $1\leq |\omega_2|\leq \cdots \leq |\omega_r|$.
\end{lemma}

\begin{proof}
    The `if' part is obvious since any successive minimum basis is orthogonal. Assume $\Lambda$ is reduced. Then we can always find a successive minimum basis $(\omega_1,...,\omega_r)$ with $\omega_1=1$ and $1\leq |\omega_2|\leq \cdots \leq|\omega_r|$. For any $2\leq j\leq r$ we have
    $$|\omega_j|_{\ti}=\inf_{x\in F_{\infty}} |\omega_j-x|=\inf_{x\in F_{\infty}}\max\{|\omega_j|,|x|\}\geq |\omega_j|,$$
    by the orthogonality property.
    Thus $|\omega_j|_{\ti}=|\omega_j|\geq 1$, \textit{i.e.}, $\omega_j\in \mathcal{F}$ for all $2\leq j\leq r$. This proves the `only if' part.
\end{proof}

The following result clarifies the relationship between $\CF$ and reduced lattices.

\begin{lemma}\label{CFimag}
    Let $z\in\Omega$. The statements (1)--(5) below are equivalent.
    \begin{enumerate}[1.]
        \item $\Lambda_z$ is reduced.
        \item $|z|_\ti \geq 1$.
        \item There exists $b\in A$ such that $z+b \in \CF$.
        \item There exists $z'\in\CF$ such that $\Lambda_z = \Lambda_{z'}$.
        \item For any $\Gamma$-representative $\tilde{z}\in\CF$ of $z$, we have $|\tilde{z}|_\ti = |z|_\ti$. 
    \end{enumerate}
    For all $z\in\Omega$, we also have
    % If the above equivalent conditions hold, then we also have
    \begin{enumerate}[6.]
        \item The covolume of $\Lambda_z$ is $D(\Lambda_z) = |z|_\ti$.
    \end{enumerate}
\end{lemma}

\begin{proof}
    (1) $\Rightarrow$ (2): Suppose $\Lambda_z$ is reduced. Then by Lemma \ref{smb} it has a successive minimum basis of the form $(1, \omega)$ with $\omega\in\CF$. We may write
    \[
    \omega = az+b, \quad\text{and}\quad z = c\omega +d, \quad\text{with}\quad a,b,c,d\in A.
    \]
    Composing these two expressions, we see that $ac=1$ and in particular $|c|=1$.

    Now 
    \[
    |z|_\ti = \inf_{x\in F_\infty}|c\omega + d + x| = 
    |c|\inf_{x\in F_\infty}|\omega + d/c + x/c| = |c||\omega|_\ti = |\omega|_\ti 
    = |\omega| \geq 1.
    \]
    In this case, we also see that $D(\Lambda_z) = |\omega| = |z|_\ti$, so (1) $\Rightarrow$ (6).

    (2) $\Rightarrow$ (3): Suppose that $|z|_\ti \geq 1$. Write $|z|_\ti = |z + x|$ with
    $x\in F_\infty$. Since $|z|\geq|z|_\ti\geq 1$, we may choose $x=b\in A$.

    Now $|z|_\ti = |z + b| = |z + b|_\ti \geq 1$ and $z+b \in\CF$.

    (3) $\Rightarrow$ (4) is obvious.

    (4) $\Rightarrow$ (5): Suppose $\Lambda_z = \Lambda_{z'}$ with $z'\in\CF$.
    Then $z = az'+b$ and $z' = cz +d$, with $a,b,c,d\in A$. Composing shows that 
    $ac=1$ and in particular $|a|=1$. 

    Now $|z|_\ti = |az'+b|_\ti = |a||z' + b/a|_\ti = |z'|_\ti$. 

    Let $\tilde{z}\in\CF$ be a $\Gamma$-representative of $z$, then $z'$ and $\tilde{z}$ are $\Gamma$-equivalent in $\CF$ and $|\tilde{z}|_\ti = |z'|_\ti = |z|_\ti$, by Proposition \ref{fundamenta domain}.

    (5) $\Rightarrow$ (2) is obvious.

    (2) $\Rightarrow$ (1): Suppose $|z|_\ti \geq 1.$ We already have $1 \in \Lambda_z$, now let $0 \neq \lambda = az + b\in\Lambda_z$. If $a=0$ then $|\lambda| \geq 1$.  Suppose $a\neq 0$, then 
    $|\lambda| = |az + b| = |a||z + b/a| \geq |z|_\ti \geq 1$. It follows that $\Lambda_z$ is reduced.

    It remains to show that (6) holds even when $\Lambda_z$ is not reduced. 
    Let 
    $\gamma=\begin{pmatrix}
    a & b \\ c & d
    \end{pmatrix}\in\GL_2(A)$  so that $\tilde{z} = \gamma(z)\in\CF$ and $\Lambda_{\tilde{z}}$ is reduced.
    Now
    \[
    D(\Lambda_z) = |cz+d|^2 \cdot D(\Lambda_{\tilde{z}}) = |cz+d|^2\cdot |\tilde{z}|_i = |cz+d|^2\cdot\frac{|z|_i}{|cz+d|^2} = |z|_i.
    \]
\end{proof}

Given $z\in\Omega$, we call $\tilde{z}$ a {\em reduced representative} of $z$ if it is $\Gamma$-equivalent to $z$ and $\Lambda_{\tilde{z}}$ is a reduced lattice.
For example, $\tilde{z}$ might be a $\Gamma$-representative of $z$ in $\CF$, in which case we also call it a {\em fundamental representative} of $z$.

% \begin{lemma}\label{covolume}
%     Let $z\in\Omega$. Then $D(\Lambda_z) = |z|_\ti$.
% \end{lemma}

% \begin{proof}
%     First, suppose that $\Lambda_z$ is reduced. Then $D(\Lambda_z) = |z|_\ti$ by Lemma \ref{CFimag}.
%     In the general case, let 
%     $\gamma=\begin{pmatrix}
%     a & b \\ c & d
%     \end{pmatrix}\in\GL_2(F_\infty)$  so that $\tilde{z} = \gamma(z)\in\CF$ is the reduced representative of $z$. 
%     Now
%     \[
%     D(\Lambda_z) = |cz+d|^2 \cdot D(\Lambda_{\tilde{z}}) = |cz+d|^2\cdot |\tilde{z}|_i = |cz+d|^2\cdot\frac{|z|_i}{|cz+d|^2} = |z|_i.
%     \]
% \end{proof}

% \begin{ex}[TO BE CHECKED]
%    The statement of Lemma \ref{CFimag}.6 may fail if $|z|_\ti < 1$.
%     Let $z=t^{-1}+t^{-3}+ut^{-6}$ with $u\in \overline{\BF_q}\smallsetminus \BF_q$, then $|z|_{\mathrm{i}}=q^{-6}$, so $\Lambda_z$ is not reduced. 
%     One can check that $\omega = t^3z-t^2-1=ut^{-3}$ is a minimal element of $\Lambda_z$.
%    Indeed, for any $a,b\in A$ with $b\neq 0$, we have $|aut^{-3}+b|\geq 1$. 
%     Thus $(\omega, 1)$ is a successive minimal basis of $\Lambda_z$ and $D(\Lambda_z)= q^{-3}> |z|_{\mathrm{i}}$.
% \end{ex}

% \begin{rem}[TO BE CHECKED]
%     The situation is different for complex elliptic curves. If $\tau\in\BH$ lies in the Poincar\'e upper half plane, then the area of a fundamental parallelogram of the lattice $\tau\BZ + \BZ\subset \BC$ is always given by $\mathrm{Im}(\tau)$, even if $\tau$ does not lie in the fundamental domain.
% \end{rem}

\subsection{Drinfeld modular forms}\label{DMF}
For a background on Drinfeld modular forms, see for example \cite{GekTur}. We will follow normalization conventions consistent with \cite{B16} and the second half of \cite{Gek88}.

To a lattice $\Lambda\subset\Cinf$ we associate the exponential function
\[
e_\Lambda(z) := z\prod_{\lambda\in\Lambda\smallsetminus\{0\}}\left(1 - \frac{z}{\lambda}\right).
\]
The function $e_\Lambda : \Cinf \rightarrow\Cinf$ is entire, surjective, $\BF_q$-linear and has simple zeros exactly at the elements of $\Lambda$. 

For every lattice $\Lambda\subset\Cinf$ there is a Drinfeld module $\phi^\Lambda$ satisfying
\[
e_\Lambda(az) = \phi_a^\Lambda(e_\Lambda(z))\quad\text{for all}\quad a\in A.
\]
In particular, following Theorem 5.4.7 of \cite{Papikian} (based on (5.4.12)-(5.4.11)-(5.4.1) in \cite{Papikian}), let 
\[
\bar\pi = \sqrt[q-1]{t-t^q}\prod_{i=1}^\infty \left(1 - \frac{t^{q^i}-t}{t^{q^{i+1}}-t}\right),
\]
then the rank 1 lattice $\Lambda = \bar\pi A$ corresponds to the Carlitz module $\rho := \phi^{\bar\pi A}$, characterised by
\[
\rho_t =  t + \tau.
\]
We also refer the reader to \cite{Gek86} (4.10) for another formula for $\bar{\pi}$.

Let $k, m\in\BZ$ and recall that a \textit{Drinfeld modular form of weight $k$ and type $m$} for the group $\Gamma=\GL_2(A)$ is a rigid analytic function $f : \Omega \rightarrow \Cinf$ satisfying
\begin{enumerate}
    \item[(A)] for all $\gamma\in\Gamma$, $f(\gamma(z)) = (cz+d)^k\det(\gamma)^{-m} f(z)$, where $(c,d)$ is the bottom row of the matrix $\gamma$; and
    \item[(B)] $f(z)$ is holomorphic at infinity.
\end{enumerate}
Condition (B) above is best understood as follows. We define a parameter at infinity 
\[
u(z) := \frac{1}{e_{\bar\pi A}(\bar\pi z)}
= \frac{1}{\bar\pi e_A(z)}.
\]
If $f$ satisfies the weak modularity condition (A), then it can be shown that $f$ has a $u$-expansion of the form
\[
f(z) = \sum_{n\in \BZ}c_n u(z)^n, \quad c_n\in\Cinf.
\]
Now condition (B) states that $c_n = 0$ for all $n < 0$.

Let $z\in\Omega$ and define the rank 2 lattice $\bar{\pi}\Lambda_z = \bar\pi(z A + A)$. This gives rise to a rank 2 Drinfeld module $\phi^z = \phi^{\bar{\pi}\Lambda_z}$ characterised by
\[
\phi^z_t = t + g(z)\tau + \Delta(z)\tau^2.
\]
It turns out that the functions $z\mapsto g(z)$ and $z\mapsto \Delta(z)$ are Drinfeld modular forms of type 0 and weights $q-1$ and $q^2-1$, respectively. Furthermore, $\Delta$ is a cusp form which vanishes nowhere on $\Omega$. 

Another important Drinfeld modular form is Gekeler's $h$-function \cite{Gek88} (we keep the traditional notation for Gekeler's $h$-function despite the conflict with our height functions, but the context will always be clear). It has weight $q+1$ and type 1 and satisfies 
\[
h(z)^{q-1} = -\Delta(z).
\]
On the other hand, the $j$-invariant of $\phi^z$, defined by $j(z) := g(z)^{q+1}/\Delta(z)$, defines a modular function
\[
j : \Omega \rightarrow \Cinf,
\]
which is modular of weight $0$ (\textit{i.e.} $\Gamma$-invariant),
but is not holomorphic at infinity.

Two points $z,z'\in\Omega$ correspond to isomorphic Drinfeld modules over $\Cinf$ if and only if $z' = \gamma(z)$ for some $\gamma\in \Gamma=\GL_2(A)$. 

If we merely have $z' = \gamma(z)$ with $\gamma\in\GL_2(F)$, then the corresponding Drinfeld modules are isogeneous. Since every isogeny $f : \phi \rightarrow\phi'$ factors as $f = \phi_a\circ f'$ with $f' : \phi\rightarrow\phi'$ cyclic, \textit{i.e.} $\ker f' \cong A/NA$ for some $a, N\in A$, we are particularly interested in cyclic isogenies.

Let $N\in A$ be a monic polynomial. 
Consider the following finite set of matrices:
$$C_N:=\left\{
\begin{pmatrix}
    a & b \\
    0 & d
\end{pmatrix}: a,b,d\in A, \;ad=N \text{ with $a,d$ monic, $|b|<|d|$ and $\gcd(a,b,d)=1$.}
\right\}$$

The number of elements of $C_N$ is 
\[
\#C_N = \sum_{d|N}\frac{|d|\cdot\varphi(e_d)}{|e_d|}=|N|\prod_{P|N}\left(1+\frac{1}{|P|}\right)
= \psi(N),
\]
where the above sum is over all monic divisors of $N$, we denote $e_d=\gcd(d,N/d)$, and the Euler totient function $\varphi$ is given by
$$\varphi(N)=|N|\prod_{P|N}\left(1-\frac{1}{|P|}\right).$$

The significance of $C_N$ is that, up to isomorphism, all targets of cyclic $N$-isogenies $f : \phi^z \rightarrow \phi^{z'}$ are are given by $z' = \gamma(z)$ for $\gamma\in C_N$. 

\begin{definition}\label{drmopoly}
    The \textit{Drinfeld modular polynomial} $\Phi_N(X,Y)\in \Cinf[X,Y]$ is the one obtained by replacing each occurrence of $j$ in the minimal polynomial of the function $z \mapsto j(N\cdot z)$ over $\Cinf(j)$ by an indeterminate $Y$. 
\end{definition}
It was shown in \cite{Bae92} that $\Phi_N(X,Y)$ has coefficients in $A$. It is also irreducible in $\Cinf[X,Y]$, symmetric in $X$ and $Y$ when $\deg N > 0$ and has degree $\psi(N)$ in each variable. Most importantly,
\[
\Phi_N\big(j(z), j(z')\big) = 0
\]
if and only if there exists a cyclic $N$-isogeny $f : \phi^z \rightarrow \phi^{z'}$.

In particular, we have
\[
\Phi_N\big(X, j(z)\big) = \prod_{\gamma\in C_N}\big(X - j(\gamma(z))\big).
\]

We will need the following result.
\begin{prop}\label{davg}
    For all $z\in\Omega$, we have
    \[
    \sum_{\gamma\in C_N}\log|\Delta(\gamma(z))| = \psi(N)\log|\Delta(z)|.
    \]
\end{prop}

This proposition is an immediate consequence of the following result, which appears in \cite[Cor 2.11]{Gek0} for the case $N$ prime, that lets us think of Gekeler's $h$ (and thus also $\Delta$) as a ``multiplicative Hecke eigenform''.
\begin{thm}\label{multhecke}
    Gekeler's $h$-function satisfies
    \[
    \prod_{\gamma\in C_N}h(\gamma(z)) = h(z)^{\psi(N)}\quad\text{for all}\quad z\in\Omega.
    \]
\end{thm}

\begin{proof}
    Denote the left hand side by $f(z) := \prod_{\delta\in C_N}h(\delta(z))$. We will show that $f(z)$ is a Drinfeld modular form of weight $\psi(N)(q+1)$ and type $\psi(N)$. As $h(z)$ is non-zero on $\Omega$, it follows that $f(z)/h(z)^{\psi(N)}$ is a modular form of weight 0, hence constant. The value of this constant is 1, which will follow from the $u$-expansion of $f(z)$.

    For a $2\times 2$ matrix $\gamma = \begin{pmatrix}
        a & b \\ c & d 
    \end{pmatrix}$, denote by
    \[
    \alpha_{k,m}(\gamma, z) := (cz + d)^k\det(\gamma)^{-m}
    \]
    the factor of automorphy, so the functional equation for a modular form of weight $k$ and type $m$ can be written $f(\gamma(z)) = \alpha_{k,m}(\gamma, z)f(z)$ for all $\gamma\in\Gamma$. A simple calculation shows that
    \[
    \alpha_{k,m}(\gamma\delta, z) = \alpha_{k,m}(\gamma, \delta(z))\cdot \alpha_{k,m}(\delta, z).
    \]

    We first show that $f(z)$ satisfies the correct functional equations. Note that by \cite[Theorem 1]{Bae92}, for every $\gamma\in\Gamma$ and $\delta\in C_N$, there exists $\gamma'\in\Gamma$ (depending on $\delta$) such that $\delta\gamma = \gamma'\delta$.

    Now let $\gamma\in\Gamma$. Then
    {\allowdisplaybreaks
\begin{align*}
        f(\gamma(z)) &= \prod_{\delta\in C_N} h(\delta\gamma(z)) 
        = \prod_{\delta\in C_N} h(\gamma'\delta(z)) \\
        &= \prod_{\delta\in C_N} \alpha_{q+1, 1}(\gamma', \delta(z)) h(\delta(z)) \\
        &= \left[\prod_{\delta\in C_N}\frac{\alpha_{q+1, 1}(\gamma'\delta, z)}{\alpha_{q+1,1}(\delta, z)}   \right] f(z) \\
        &= \left[\prod_{\delta\in C_N}\frac{\alpha_{q+1, 1}(\delta\gamma, z)}{\alpha_{q+1,1}(\delta, z)}   \right] f(z) \\
        &= \left[\prod_{\delta\in C_N}\frac{\alpha_{q+1, 1}(\delta, \gamma(z))}{\alpha_{q+1,1}(\delta, z)}\alpha_{q+1,1}(\gamma, z)   \right] f(z) \\
        &= \alpha_{q+1, 1}(\gamma, z)^{\#C_N}f(z) = 
        \alpha_{\psi(N)(q+1), \psi(N)}(\gamma, z)f(z),
    \end{align*}}
as required, since for $\delta\in C_N$ we see that $\alpha_{k,m}(\delta,z) = d^kN^{-m}$ is independent of $z$.

Next, we compute the first term in the $u$-expansion of $f(z)$, showing that it is indeed holomorphic at infinity and thus a modular form. We will actually work with the following parameter, associated to the full modular group $\Gamma(N) = \ker\big(\Gamma \rightarrow \GL_2(A/NA)\big)$:
\[
u_N(z) := \frac{1}{e_{\bar\pi A}(\bar\pi z/N)}.
\]
It is related to $u(z) = u_1(z)$ as follows. Note that $\rho_N(X) = X^{|N|} + \cdots + NX$, as $N$ is monic.
{\allowdisplaybreaks
\begin{align*}
    u_1(z) &= \frac{1}{e_{\bar\pi A}(\bar\pi z)} = \frac{1}{\rho_N(e_{\bar\pi A}(\bar\pi z/N))} \\
    &= \frac{1}{e_{\bar\pi A}(\bar\pi z/N)^{|N|} + \cdots + Ne_{\bar\pi A}(\bar\pi z/N)} \\
    &= \frac{1}{u_N(z)^{-|N|} + \cdots + Nu_N(z)^{-1}}\\
    &= u_N(z)^{|N|}\big[1 + O(u_N(z))\big].
\end{align*}}

The $u$-expansion of $h(z)$ starts with
\[
h(z) = -u(z) + O(u(z)^2) = -u_N(z)^{|N|} + O(u_N(z)^{|N|+1}).
\]

Now, for $\delta = \begin{pmatrix}
    a & b \\ 0 & d
\end{pmatrix}\in C_N$,
we compute 
\begin{align*}
    u_1\left(\frac{az + b}{d}\right) &= u_1\left(\frac{a^2z + ab}{N}\right) 
    = u_N(a^2z + ab) \\
    &= \frac{1}{e_{\bar\pi A}(\bar\pi(a^2z+ab)/N)} 
    = \frac{1}{e_{\bar\pi A}(\bar\pi a^2z/N) + e_{\bar\pi A}(\bar\pi ab/N)}\\
    &= \frac{1}{\rho_{a^2}(e_{\bar\pi A}(\bar\pi z/N)) + \rho_{ab}(e_{\bar\pi A}(\bar\pi/N))}\\
    &= \frac{1}{u_N(z)^{-|a^2|} + \cdots + a^2 u_N(z)^{-1} + \rho_{ab}(e_{\bar\pi A}(\bar\pi/N))}\\
    &= u_N(z)^{|a^2|} + O(u_N(z))^{|a^2|+1}.
\end{align*}

Thus, $f(z)$ has a $u_N$-expansion starting with
\begin{align*}
    f(z) &= \prod_{\delta\in C_N}h(\delta(z)) \\
    &= \prod_{\delta\in C_N} \left[ -u_1((az+b)/d) + O(u_1((az+b)/d)^2) \right]\\
    &= \prod_{\delta\in C_N} \left[ -u_N(z)^{|a^2|} + O(u_N(z))^{|a^2|+1} \right]\\
    &= (-1)^{\psi(N)}u_N(z)^{|N|\psi(N)} + O\big(u_N(z)^{|N|\psi(N)+1}\big),
\end{align*}
where we verify that
\[
\sum_{\delta\in C_N}|a|^2 = |N|\sum_{\delta\in C_N}\frac{|a|}{|d|} 
= |N|\sum_{d|N}\frac{\varphi(e_d)}{|e_d|}|a| = |N|\psi(N). 
\]
Since the $u_N$-expansion of $f$ has no polar terms, the same holds for its $u$-expansion, hence it is a modular form for $\Gamma$. We also see $f(z)$ and $h(z)^{\psi(N)}$ have the same leading term in their $u_N$-expansions. The result follows.
\end{proof}

\subsection{Bruhat-Tits trees}\label{sec:Tree}

Given a non-trivial vector space $V$ over a local field, one can associate a simplicial complex of dimension $\dim(V)-1$ to $\PGL(V)$, which is called the \textit{Bruhat-Tits building} for $\PGL(V)$. It is analogous to the symmetric space for real Lie groups. In this article, we consider the vector space $V=F_\infty\oplus F_\infty$ so that the associated Bruhat-Tits building becomes a tree, which we shall call the \textit{Bruhat-Tits tree}. Let us give a precise definition of the Bruhat-Tits tree.

In this section, it is convenient to denote by $\pi = t^{-1}$ the uniformiser of $F_\infty$.

\begin{definition}
    The \textit{Bruhat-Tits tree} $\Tcal$ for $V:=F_\infty\oplus F_\infty$ is a $(q+1)$-regular tree with vertices $v_i$ and edges $e_{i,j}$. 
    Each vertex $v_i$ of $\Tcal$ represents a homothety class of (rank two) $\Ocal_\infty$-lattices in $V$, and an edge $e_{i,j}$ connects vertices $v_i$ and $v_j$ if and only if there exist $\Ocal_\infty$-lattices $L_i$ and $L_j$ in the classes of $v_i$ and $v_j$, respectively, such that $\pi L_j\subsetneq L_i\subsetneq L_j$.
\end{definition}

Let $g\in \GL_2(F_\infty)$ and $L$ be any $\Ocal_\infty$-lattice in $V$. We fix a left action of $\Gamma$ on the set of $\Ocal_\infty$-lattices given by $g\cdot L:=Lg^{-1}$. We denote by $X(\Tcal)$ (resp. $Y(\Tcal)$) the set of vertices (resp. directed edges) of $\Tcal$. They are identified by
\begin{align*}
    \GL_2(F_\infty)/F_\infty\cdot \GL_2(\Ocal_\infty) & \xrightarrow{\sim} X(\Tcal),\\
    g &\mapsto [g\cdot L_0];\\
    \GL_2(F_\infty)/F_\infty\cdot \Ifra & \xrightarrow{\sim} Y(\Tcal),\\
    g &\mapsto g\cdot e_0,
\end{align*}
where $L_i$ is the $\Ocal_\infty$-lattice $\pi^{-i}\Ocal_\infty\oplus \Ocal_\infty$ and $e_0$ is the edge originating from $L_0$ to $L_{-1}$, and $\Ifra$ is the \textit{Iwahori group} given by
$$\Ifra:=\left\{\begin{pmatrix}
    a & b\\
    c & d
\end{pmatrix}\in \GL_2(\Ocal_\infty): c\equiv 0\ (\text{mod }\pi)\right\}.$$
Following \cite{Gek0} page 184, one can see that $X(\Tcal)$ can be represented by the set of matrices:
$$\left\{\begin{pmatrix}
    \pi^k & u\\
    0 & 1
\end{pmatrix}:k\in \BZ,\ u\in F_\infty/\pi^k\Ocal_\infty\right\}.$$
As the action by $\begin{pmatrix}
    0 & \pi^{-1}\\
    1 & 0
\end{pmatrix}$ reverses the direction of each $e\in Y(\Tcal)$, the set of matrix representatives of $Y(\Tcal)$ is given by (see also \cite{Gek0} page 184)
$$\left\{\begin{pmatrix}
    \pi^k & u\\
    0 & 1
\end{pmatrix}:k\in \BZ,\ u\in F_\infty/\pi^k\Ocal_\infty\right\}\sqcup  \left\{\begin{pmatrix}
    \pi^k & u\\
    0 & 1
\end{pmatrix}:k\in \BZ,\ u\in F_\infty/\pi^k\Ocal_\infty\right\}\cdot\begin{pmatrix}
    0 & \pi^{-1}\\
    1 & 0
\end{pmatrix}.$$
Consequently, each vertex can be parametrised by $v(k,u)$ with $k\in \BZ$ and $u\in F_\infty/\pi^k\Ocal_\infty$. In particular, the vertex $v_k:=[L_k]$ is represented by $v(k,0)$. These vertices generate the \textit{principal axis} $A(0,\infty)$ of $\Tcal$. It has two ends: $h_0$ given by $v_k$ with $k\in \BZ^{\geq 0}$ and $h_\infty$ given by $v_k$ with $k\in \BZ^{\leq 0}$. It is well-known that $h_\infty$ is a fundamental domain for $\Gamma$ in the following sense (cf. \cite[Theorem 6.4]{GekTur}):

\begin{thm}
    Each vertex of $\Tcal$ is $\Gamma$-equivalent to precisely one of the $v_k$ $(k\leq 0)$. Similarly, each directed edge is equivalent to precisely one of the $e_k$ or $\overline{e_k}$ $(k\leq 0)$, where $e_k$ is the edge originating from $v_k$ to $v_{k-1}$ and $\overline{e_k}$ is the reverse of $e_k$.
\end{thm}

However, given a vertex $v(k,u)$ of $\Tcal$, it is not immediately clear which vertex $v_k\in h_\infty$ is $\Gamma$-equivalent to it. We partially solve this problem in Theorem \ref{fdmr} below.
The following Proposition \ref{avdc} and Theorem \ref{fdmr} are not new to the experts and similar calculations can be seen in \cite[\S6.5]{Gek95}. We include the proofs here since we could not find a reference to the precise statements we need.

Let $v(k,u)\in X(\Tcal)$ be a vertex. As $u\in F_\infty/\pi^k\Ocal_\infty$, we can write $u=u_0+u_1$ with $u_0\in \BF_q[t]$ and $v_\infty(u_1)>0$. Suppose $u_1\neq 0$ modulo $\pi^k\Ocal_\infty$, we can then express $u$ uniquely as 
\begin{equation}\label{expu}
    u=u_0+u_1 \text{ with } u_0\in \BF_q[t] \text{ and } u_1=\sum_{i=1}^ra_i\pi^i \text{ such that } r<k, a_i\in \BF_q\text{ and } a_r\neq 0.
\end{equation}

\begin{prop}\label{avdc}
    Let $v(k,u)\in X(\Tcal)$ be a vertex in the Bruhat-Tits tree $\Tcal$. Suppose $u$ is expressed as in (\ref{expu}). Then the vertices adjacent to $v(k,u)$ are those $v(k',w)\in X(\Tcal)$ such that either $k'=k+1$ and $w=u+\xi\pi^k$ with $\xi\in \BF_q$, or $k'=k-1$ and $w=u\text{ mod }\pi^{k-1}\Ocal_\infty$. 
\end{prop}

\begin{proof}
    Let $L(k,u)$ be the $\Ocal_\infty$-lattice representing $v(k,u)$ given by $(\pi^{-k},-u\pi^{-k})\Ocal_\infty\oplus (0,1)\Ocal_\infty$, similarly for $L(k',w)$. Let $w=u+\xi\pi^k$ with any $\xi\in \BF_q$. One can check that 
    $$\pi(\pi^{-(k+1)},-w\pi^{-(k+1)})=(\pi^{-k},-u\pi^{-k})-\xi(0,1),$$
    from which one can deduce $\pi L(k+1,w)\subsetneq L(k,u)$. Also, one can check that 
    $$(\pi^{-k+1},-u\pi^{-k+1})=\pi(\pi^{-k},-w\pi^{-k})+\xi\pi(0,1),$$
    from which one can again deduce that $\pi L(k,u)\subsetneq \pi L(k+1,w)$. Therefore, we have obtained 
    $$\pi L(k,u)\subsetneq \pi L(k+1,w)\subsetneq L(k,u),$$
    which is equivalent to saying $v(k,u)$ is adjacent to $v(k+1,w)$. 

    The same calculation also gives
    \[
    \pi L(k-1, w) \subsetneq \pi L(k, u) \subsetneq L(k-1, w)
    \]
    for $w = u \bmod \pi^{k-1}\Ocal_\infty$.
    This completes the proof.
\end{proof}

\begin{thm}\label{fdmr}
    Let $v(k,u)\in X(\Tcal)$ be a vertex, and $u=u_0+u_1$ be expressed as in (\ref{expu}) so that $\deg_{\pi}(u_1)=r$. Assume $u_1=0$ if $u=u_0\in \BF_q[t]$. Let $k'\leq 0$ be an integer such that $v_{k'}$ is the vertex on $h_\infty$ that is $\Gamma$-equivalent to $v(k,u)$. 
    \begin{enumerate}[1.]
        \item If $u=0$, then $k'=-|k|$.
        \item If $\deg_{\pi}(u_1)\geq k$ (\textit{i.e.} $u_1=0$), then $k'=-|k|$.
        \item If $0<r<k$ and $k\geq 2r$, then $k'=2r-k$.
        \item If $0<r<k$ and $k< 2r$, then $k'\geq k-2r$.
    \end{enumerate}
\end{thm}
% [CHECK: the referee claims this is in \cite{Gek0}. I didn't find it. Maybe cite another paper of Gekeler, for instance Improper Eisenstein paragraph 6.5...]

\begin{proof}
    The proof is constructive. Suppose $\begin{pmatrix}
        \pi^k & u\\
        0 & 1
    \end{pmatrix}$ corresponds to $v(k,u)$.

    If $u=0$, we only need to consider the case when $k>0$. Indeed, we have 
    $$\begin{pmatrix}
        \pi^{-k} & 0\\
        0 & 1
    \end{pmatrix}= 
    \begin{pmatrix}
        0 & 1\\
        1 & 0
    \end{pmatrix}
    \begin{pmatrix}
        \pi^{k} & 0\\
        0 & 1
    \end{pmatrix}
    \begin{pmatrix}
        0 & 1\\
        1 & 0
    \end{pmatrix}\pi^{-k}.$$
    This shows that $v_k$ is translated to $v_{-k}\in h_\infty$, proving the first statement. 

    If $u_1=0$, then 
    $$\begin{pmatrix}
        \pi^{k} & u\\
        0 & 1
    \end{pmatrix}= 
    \begin{pmatrix}
        u_0 & 1\\
        1 & 0
    \end{pmatrix}
    \begin{pmatrix}
        \pi^{-k} & 0\\
        0 & 1
    \end{pmatrix}
    \begin{pmatrix}
        0 & 1\\
        1 & 0
    \end{pmatrix}\pi^{k}.$$
    Therefore, combining with the first statement we prove the second statement.

    If $0<r<k$ and $k\geq 2r$, then there exist $c,d\in \BF_q[t]$ with $\deg_t(c)=r$ such that $dt^r+cu_1t^r=1$. We shall prove the existence after looking at the decomposition:
    $$\begin{pmatrix}
        \pi^{k} & u\\
        0 & 1
    \end{pmatrix}= 
    \begin{pmatrix}
        d-cu_0 & t^ru\\
        -c & t^r
    \end{pmatrix}
    \begin{pmatrix}
        \pi^{k-2r} & 0\\
        0 & 1
    \end{pmatrix}
    \begin{pmatrix}
        1 & 0\\
        c\pi^{k-r} & 1
    \end{pmatrix}\pi^{r}.$$
    Under our assumption, $c\pi^{k-r}\in \Ocal_\infty \text{ and } t^ru\in \BF_q[t]$, and one can verify that $\begin{pmatrix}
        d-cu_0 & t^ru\\
        -c & t^r
    \end{pmatrix}\in \Gamma.$ Now we show the existence of such $c,d$. Indeed, it is equivalent to finding $c,d\in \BF_q[t]$ such that $d+cu_1=\pi^r$. We may take $c=c_0+c_1t+\cdots+t^r$. Then finding solutions to the equations
    \begin{align*}
        c_0a_1+c_1a_2+\cdots+c_{r-1}a_r &=0,\\
        c_0a_2+c_1a_3+\cdots+c_{r-2}a_r &=0,\\
        \vdots\\
        c_0a_{r-1}+c_1a_r &=0,\\
        c_0a_r&=1\
    \end{align*}
    gives the existence of such $c,d$. Recursively one can find solutions to the above equations since $a_r\neq 0$, hence we have proved the existence of the desired $c,d\in \BF_q[t]$. Again combining with the first statement, we complete the proof of the third statement.

    Assume $0<r<k$ and $k<2r$. According to Proposition \ref{avdc}, $v(k,u)$ is adjacent to $v(k+1,u)$ and $v(k+1,u)$ is adjacent to $v(k+2,u)$. Continuing this process, we get a path from $v(k,u)$ to $v(2r,u)$ of length $2r-k$. By the third statement, $v(2r,u)$ is $\Gamma$-equivalent to $v_0$. Therefore, $v(k,u)$ is $\Gamma$-equivalent to some $v_{k'}$ with $k'\geq k-2r$.
\end{proof}

\subsection{The Building map}\label{sec:BuildingMap}

We denote by $\Tcal(\BR)$ the {\em realization} of $\Tcal$, where each edge $e$ of $\Tcal$ is identified with a copy of the unit interval $[0,1]$ with $0$ corresponding to the source and $1$ with the target vertex of that edge. We write $e(\BR)\cong [0,1]$ for this ``filled-in'' edge.
By a theorem of Goldmann-Iwahori, the points of $\Tcal(\BR)$ are in bijection with similarity classes of norms on $V$, see for example \cite[Theorem 1.4.3]{GeRe}.

\begin{definition}
    The \textit{building map} $\lambda$ on $\Omega$ is given by 
    \begin{align*}
        \lambda: \Omega &\rightarrow \Tcal(\BR) %\\
        % z &\mapsto |\cdot|_z: (x,y)\mapsto |xz+y|
    \end{align*}
    where $\lambda(z)\in\Tcal(\BR)$ is the point associated to the norm $\|\cdot\|_z : (x,y) \mapsto |xz+y|$ on $V$.
\end{definition}

We collect some useful results, which can be found in \cite{Gek0} and \cite{GeRe}.

\begin{prop}\label{impz}Let $z\in \Omega$ and $g\in \GL_2(F_\infty)$.
    \begin{enumerate}[1.]
        % \item $|g\cdot z|_{\ti}=|cz+d|^{-2}|\det(g)||z|_{\ti}$.
        \item \cite[Lemma 1.8]{Gek0} If $\lambda(z)\in X(\Tcal)$ and $\lambda(z)$ is represented by the matrix $\begin{pmatrix}
    \pi^k & u\\
    0 & 1
    \end{pmatrix}$, then
    $$-\log|z|_{\mathrm{i}}=k,\quad \text{and}\quad -\log|z|=
    \begin{cases}
        v_\infty(u),\ &\text{if } u\neq 0 \ (\text{mod}\ \pi^k\mathcal{O}_\infty);\\
        k, &\text{if } u=0\ (\text{mod}\ \pi^k\mathcal{O}_\infty).
    \end{cases}$$
    \item \cite[Cor. 1.9]{Gek0} We have $|z|=|z|_{\mathrm{i}}$ if and only if $\lambda(z)\in A(0,\infty)$.
    \item \cite[Prop. 1.5.3]{GeRe} The building map is invariant under the action of $\GL_2(F_\infty)$. This means that for any $g\in \GL_2(F_\infty)$, the point $\lambda(g\cdot z)\in \Tcal$ corresponds to the similarity class of the norm $\|\cdot\|_{g\cdot z}$, defined by 
    $$\|\alpha\|_{g\cdot z} = \|\alpha g\|_z, \quad \forall \alpha\in V.$$
    If $\lambda(z)=[L]\in X(\Tcal)$ with an $\Ocal_\infty$-lattice representative $L$, then $\lambda(g\cdot z) = [g\cdot L]$.
    % $$\lambda(g\cdot z)(\alpha) \text{ is similar to the norm }
    % g\|\cdot\|_z(\alpha):=\|\alpha g\|_z,\ \forall \alpha\in V, \forall g\in \GL_2(F_\infty);$$
    % and if $\lambda(z)=[L]\in X(\Tcal)$ with an $\Ocal_\infty$-lattice representative $L$, then
    % $$\lambda(g\cdot z) \sim g\|\cdot\|_z \text{ corresponds to the lattice class } [g\cdot L].$$
    \end{enumerate}
\end{prop}

% We collect some useful results, which can be found in \cite{Gek0},\cite{GeRe}.

% \begin{prop}\label{impz}Let $z\in \Omega$ and $g\in \GL_2(F_\infty)$.
%     \begin{enumerate}[1.]
%         % \item $|g\cdot z|_{\ti}=|cz+d|^{-2}|\det(g)||z|_{\ti}$.
%         \item (see \cite[Lemma 1.8]{Gek0}) If $\lambda(z)\in X(\Tcal)$ and $\lambda(z)$ is represented by the matrix $\begin{pmatrix}
%     \pi^k & u\\
%     0 & 1
%     \end{pmatrix}$, then
%     $$-\log|z|_{\mathrm{i}}=k,\quad \text{and}\quad -\log|z|=
%     \begin{cases}
%         v_\infty(u),\ &\text{if } u\neq 0 \ (\text{mod}\ \pi^k\mathcal{O}_\infty);\\
%         k, &\text{if } u=0\ (\text{mod}\ \pi^k\mathcal{O}_\infty).
%     \end{cases}$$
%     \item (see \cite[Corollary 1.9]{Gek0}) We have $|z|=|z|_{\mathrm{i}}$ if and only if $\lambda(z)\in A(0,\infty)$.
%     \item (see \cite[Proposition 1.5.3]{GeRe}) The building map is invariant under the action of $\GL_2(F_\infty)$. This means 
%     $$\lambda(g\cdot z)(\alpha)=g\|\cdot\|_z(\alpha):=\|\alpha g\|_z,\ \forall \alpha\in V, \forall g\in \GL_2(F_\infty);$$
%     and if $\lambda(z)=[L]\in X(\Tcal)$ with an $\Ocal_\infty$-lattice representative $L$, then
%     $$\lambda(g\cdot z)=g\|\cdot\|_z=[g\cdot L].$$
%     [CHECK notation for norm classes.]
%     \end{enumerate}
% \end{prop}

\begin{rem}\label{exot}
    Since $\log|z|_{\ti}=-k$, on $h_\infty$ the imaginary part $|z|_{\ti}$ will increase by one if $\lambda(z)$ moves one step in the direction of the end $h_\infty$ of $\Tcal$. Therefore, the imaginary part tends to infinity if one ``moves towards infinity''. This matches our intuition from the classical case.
\end{rem}

\begin{lemma}\label{impi}
    Let $z\in \Omega$ such that $\lambda(z)\in e(\BR)$, where $e\in Y(\Tcal)$ is a directed edge with one vertex $v(k,u)$ and another vertex $v(k-1,u)$. Then $q^{-k}\leq |z|_{\mathrm{i}}\leq q^{-k+1}$. If $\lambda(z)$ lies in the interior of $e(\BR)$, then these inequalities are strict.
\end{lemma}

\begin{proof}
    This is \cite[Lemma 1.7]{Gek0}.
    It is clear that $e$ is given by the matrix $g(k,u):=\begin{pmatrix}
        \pi^k & u\\
        0 & 1
    \end{pmatrix}$. This means $e$ is the translation of $e_0$ by $g(k,u)$. Let $\tilde{z}\in\lambda^{-1}(e_0(\BR))$ be a point for which $z=g(k,u)\cdot \tilde{z}$. Thus we have
    $$|z|_{\ti}=|\det(g(k,u))|\cdot|\tilde{z}|_{\ti}=q^{-k}|\tilde{z}|_{\ti}.$$
    Since $1\leq |\tilde{z}|_{\ti}\leq q$, our assertion follows immediately.

    Finally, it follows from \cite[(1.6)]{Gek0} that $1 < |\tilde{z}|_{\ti} < q$ if $\tilde{z}$ lies in the interior of $e_0(\BR)$.
\end{proof}

\subsection{Farey sequences}\label{Farey}

Here we record Hsia's results on Farey sequences, which will be used in the sequel. Let $M\geq 1$ be an integer. The \textit{Farey sequence of order $M$} is given by 

$$F_M:=\left\{\frac{h}{f}: h,f\in \BF_q[t] \text{ with $f$ monic, } \deg(h)<\deg(f)\leq M \text{ and } \gcd(h,f)=1 \right\}$$

Let $D_M(\frac{h}{f})$ denote the set 
$$D_M\left(\frac{h}{f}\right):=\left\{\zeta\in F_{\infty}: |\zeta|<1 \text{ and } \left|\zeta-\frac{h}{f}\right|\leq \frac{1}{|f|q^{M+1}} \right\}.$$

Let $I$ be the open unit ball $\{\zeta\in F_{\infty}: |\zeta|<1\}$. Then we have

\begin{lemma}
    $I=\coprod\limits_{\frac{h}{f}\in F_M}D_M(\frac{h}{f}).$
\end{lemma}

\begin{proof}
    See \cite[Proposition 3.4]{hsia}.
\end{proof}

Let us now give a technical lemma useful in the proof of the main theorem.

\begin{lemma}\label{countD_M}
    Let $h,r \in\BF_q[t]$, let $d,e,f \in\BF_q[t]-\{0\}$, let $M\geq 1$, with $\deg(h)<\deg(f)\leq M$, and $\gcd(h,f)=1$, and such that $\frac{|d|}{|e||f|q^M}\geq1$. Then the number of elements $b\in\BF_q[t]$ satisfying
    \[
    \frac{b}{d} \in D_M\left(\frac{h}{f}\right)\quad\text{and}\quad b\equiv r \bmod e
    \]
    equals 
    \[
    \frac{|d|}{|e||f|q^M}.
    \]
    In particular, the number of $b\in\BF_q[t]$ with $\frac{b}{d}\in D_M\left(\frac{h}{f}\right)$ and $\gcd(b,e)=1$ equals
    \[
    \frac{\varphi(e)|d|}{|e||f|q^M}.
    \]
\end{lemma}

\begin{proof}
    The second assertion clearly follows from the first.

    Write $b = r + b'e$. We find that
    \[
    \frac{b}{d} \in D_M\left(\frac{h}{f}\right) \Longleftrightarrow 
    \left|\frac{b}{d} - \frac{h}{f}\right| \leq \frac{1}{|f|q^{M+1}}
    \Longleftrightarrow
    \left|b' - \frac{\frac{h}{f}d-r}{e}\right| \leq \frac{|d|}{|e||f|q^{M+1}}.
    \]
    Now the result follows from the fact that, for any $x\in F_\infty$, the number of $b'\in\BF_q[t]$ satisfying $|b' - x| \leq q^D$ is $q^{D+1}$ for $D\geq 0$.
\end{proof}

\subsection{Size of representatives}\label{Funda}

We prove a useful estimate on representatives in the Drinfeld domain, which will be important in the proof of Theorem \ref{thm:main}.

Let $z \in\CF$ and $\gamma=\begin{pmatrix}
    a & b\\
    0 & d
\end{pmatrix}\in C_N$.
To prove our main results, we need to estimate the imaginary part of a fundamental representative $\widetilde{z_\gamma} \in\CF$ of $z_\gamma := \gamma(z)$. When $|d|$ is large, this seems difficult. Instead, following Hsia \cite{hsia}, we construct an explicit ``Farey representative'' $\widehat{z_\gamma}$, which is sufficiently close to $\widetilde{z_\gamma}$ to be useful.

Suppose that $|d|>\sqrt{|N||z|_\ti}$. Set $M:=\left\lceil \log\frac{|d|}{\sqrt{|N||z|_\ti}}\right\rceil\geq 1.$

Take $\delta=\begin{pmatrix}
    \alpha & \beta \\
    f & -h
\end{pmatrix}\in \Gamma$ such that $\frac{h}{f}\in F_M$ and $\frac{b}{d}\in D_M(\frac{h}{f})$. Set $\widehat{z_{\gamma}}:=\delta \cdot z_{\gamma}$. Then using $|z|=|z|_\ti$ we compute

$$|\widehat{z_{\gamma}}|_\ti =\frac{1}{|fz_{\gamma}-h|^2}\left|\frac{a}{d}\right||z|_\ti=\frac{|d|^2}{|f|^2|N||z|_\ti}\cdot\frac{1}{\left|1+\frac{b/d-h/f}{Nz/d^2}\right|^2}.$$

\begin{lemma}\label{zhat}
    With notations as above, the following properties hold.
    \begin{enumerate}[1.]
        \item $|\widehat{z_{\gamma}}|_{\mathrm{i}}\geq \frac{1}{q^2}.$
        \item $\log|\widehat{z_{\gamma}}|_{\mathrm{i}}\leq \log \frac{|d|^2}{|N||z|_{\mathrm{i}} |f|^2}$.
        \item Let $\widetilde{z_\gamma}\in \mathcal{F}$ be a fundamental representative of $z_\gamma$. Then $\log|\widetilde{z_\gamma}|_{\mathrm{i}}\leq \log|\widehat{z_\gamma}|_{\mathrm{i}}+4.$
    \end{enumerate}
\end{lemma}

\begin{proof}
    By the definition of $M$ and the choice of $f$ we have 
    \begin{equation}\label{ubf}
        \log|f|\leq 1+\log\frac{|d|}{\sqrt{|N||z|_\ti}}, \quad\text{hence}\;\; 
        \frac{|d|^2}{|f|^2|N||z|_\ti}\geq \frac{1}{q^2}.
    \end{equation}
    To proceed, note that
    \[
    \left|1+\frac{b/d-h/f}{Nz/d^2}\right|
    = \frac{1}{|z|}\left|z+\frac{b/d-h/f}{N/d^2}\right|
    \geq \frac{1}{|z|}|z|_\ti = 1,
    \]
    since $z\in\CF$. Part 2 now follows.

    To show Part 1, either 
    \[
    \left|1+\frac{b/d-h/f}{Nz/d^2}\right| = 1,
    \]
    in which case 1 follows from (\ref{ubf}), or 
    \[
    \left|1+\frac{b/d-h/f}{Nz/d^2}\right| > 1.
    \]
    In the latter case, by  $b/d\in D_M(h/f)$, we have
    \[
    \left|\frac{b}{d}-\frac{h}{f}\right|\leq \frac{1}{|f|q^{M+1}}\implies 1<\left|\frac{b/d-h/f}{Nz/d^2}\right|\leq \frac{1}{q|f|}\frac{|d|}{\sqrt{|N||z|_\ti}}.
    \]
    Therefore
    \[
    |\widehat{z_{\gamma}}|_\ti=\frac{|d|^2}{|f|^2|N||z|_\ti}\cdot\frac{1}{\left|\frac{b/d-h/f}{Nz/d^2}\right|^2}\geq \frac{1}{q^2}.
    \]
    
    We now prove Part 3. It is here that we need our preparations from Sections \ref{sec:Tree} and \ref{sec:BuildingMap} involving the building map $\lambda:\Omega\rightarrow \Tcal(\BR)$. Suppose $e\in Y(\Tcal)$ is the edge such that $\lambda(\widehat{z_\gamma})\in e(\BR)$. Let $v(k,u)\in X(\Tcal)$ be one vertex of $e$ with the other vertex given by $v(k+1,w)$ with $w=u+\xi\pi^k$ modulo $\pi^{k+1}\Ocal_\infty$ (cf. Proposition \ref{avdc}).
    
    Thus, Lemma \ref{impi} implies $\log |\widehat{z_{\gamma}}|_\ti\leq -k$. By our first statement, this implies $k\leq 2$. 
    Let $v_{k'}\in h_\infty$ be the fundamental representative of $v(k,u)$, and write $u=u_0+u_1$ as in (\ref{expu}). 
    % In the case when $k<2$ we always have $u_1=0$ since $v_\infty(u_1)>0$ and $0<v_\infty(u_1)<k$ if $u_1\neq 0$; hence $k'=-|k|$ for $k<2$ or $k=2$ and $u_1=0$ by Theorem \ref{fdmr}. If $k=2$ and $u_1\neq 0$, then $r=\deg_{\pi}(u_1)<k=2$, which implies $r=1$. In this case, $k'=2r-k=0$ by Theorem \ref{fdmr}. 

    % The fundamental representative $\lambda(\tilde{z}_\gamma)$ of 

    To complete the proof, we distinguish three cases:

    \paragraph{\bf Case 1:} $k=2, \; u_1 = 0$.  Since $\log |\hat{z}_\gamma|_{\mathrm{i}} \leq -k = -2$ from Lemma \ref{impi} and $\log |\hat{z}_\gamma|_{\mathrm{i}} \geq -2$ from Part 1, we must have $\log |\hat{z}_\gamma|_{\mathrm{i}} = -2$ and since the inequalities in Lemma \ref{impi} are strict in the interior of the edge, it follows that 
    $\lambda(\hat{z}_\gamma) \in v(k,u)$ lies on this vertex. Its representative vertex is $v_{k'}$ with $k'=-2$ by Theorem \ref{fdmr}, so we find that
    \[
    \log|\tilde{z}_\gamma|_{\mathrm{i}} = -k' = 2 = \log |\hat{z}_\gamma|_{\mathrm{i}} + 4.
    \]

    \paragraph{\bf Case 2:} $k=2, \; u_1\neq 0$. As above, we find that $\lambda(\hat{z}_\gamma) \in v(k,u)$ lies on the vertex. Since $0 < r=v(u_1) < k =2$, we must have $r=1$ and Theorem \ref{fdmr} gives $k'=0$, so we get 
    \[
    \log|\tilde{z}_\gamma|_{\mathrm{i}} = -k' = 0 = \log |\hat{z}_\gamma|_{\mathrm{i}} + 2.
    \]

    \paragraph{\bf Case 3:} $k\leq 1$. In this case, we must have $u=0$ or $u_1=0$, since otherwise $0 < r=v(u_1) < k \leq 1$ is impossible. Thus Theorem \ref{fdmr} gives $k' = -|k|$. This means that $\lambda(\tilde{z}_\gamma)$ lies on an edge in $h_\infty$ with one vertex $v_{k'}$; the other vertex is $v_{k'\pm 1}$. Lemma \ref{impi} gives
    \[
    -k-1 \leq \log|\hat{z}_\gamma|_{\mathrm{i}} \quad\text{and}\quad \log|\tilde{z}_\gamma|_{\mathrm{i}} \leq -k' + 1 = |k| + 1.
    \]
    The result now follows since $|k|+1 \leq (-k-1) + 4 \leq \log|\hat{z}_\gamma|_{\mathrm{i}} + 4$ for all $k \leq 1$.
    
    % We complete our proof by showing a case for example. The remaining cases can be done in the same way 
    
    % [CHECK: the referee asks us to do it]. 
    % Consider the case when $k=2$ and $u_1=0$. Using Theorem \ref{fdmr} again we see $v(k+1,w)$ has the fundamental representative $v_{-1}$. Proposition \ref{impz} and Lemma \ref{impi} then give us $\log|\widetilde{z_\gamma}|_{\mathrm{i}}\leq \log|\widehat{z_\gamma}|_{\mathrm{i}}+4.$
\end{proof}

\section{Proof of Theorem \ref{thm:main}}\label{proof of main}

Let us start with some useful lemmas. The first tells us that, in the non-Archimedean setting, the height of a (monic, univariate) polynomial equals its Mahler measure.

\begin{lemma}\label{mahler}\cite[Lemma 2.2]{hsia}
    Given any $z\in \Omega$, we have
    $$h(\Phi_N(X,j(z)))=\sum_{\gamma\in C_N}\log\max \{1,|j(z_{\gamma})|\}.$$
\end{lemma}

\begin{lemma}\label{imagsum}
    Let $N\in \BF_q[t]$ be monic, with $\deg N>0$. Let $\gamma=\begin{pmatrix}
    a_\gamma & b_{\gamma}\\
    0 & d_{\gamma}
\end{pmatrix}\in C_N$.
    \begin{enumerate}[1.]
    \item We have
    \[
    \sum_{\gamma\in C_N}\log\frac{|d_\gamma|}{|a_\gamma|} = \psi(N)\big(\deg N - 2\lambda_N\big),
    \]
    where
    \[
    \lambda_N := \sum_{P^n \| N}\frac{|P|^n-1}{|P|^{n-1}(|P|^2-1)}\deg P.
    \]
    
    \item In particular, we have
    \[
    -\sum_{\gamma\in C_N} \log |z_\gamma|_{\mathrm{i}} = \psi(N)\big(\deg N - 2\lambda_N - \log |z|_{\mathrm{i}}\big).
    \]
    
    \item  We also have
    \[
    -\sum_{\gamma\in C_N}\log D(\Lambda_{z_\gamma}) = 
    \psi(N)\big(\deg N - 2\lambda_N - \log D(\Lambda_z)\big).
    \]
    \end{enumerate}
\end{lemma}

\begin{proof}
    Let $\displaystyle{S_N:=\sum_{\gamma\in C_N}\log\frac{|d_\gamma|}{|a_\gamma|}}$. Suppose $M,N\in \BF_q[t]$ are monic and coprime to each other. Then we have
    {\allowdisplaybreaks
    \begin{align*}
        S_{MN} & = \sum_{\gamma\in C_{MN}}\log\frac{|MN|}{|a_\gamma|^2}\\
        & = \sum_{\gamma\in C_{MN}}\log\frac{|MN|}{|a_{\gamma,M}|^2|a_{\gamma,N}|^2},\ \text{where $a_{\gamma,M}|M$ and $a_{\gamma,N}|N$ are coprime, and $a_\gamma=a_{\gamma,M}a_{\gamma,N}$;}\\
        &= \sum_{\gamma\in C_{MN}}\log \frac{|M|}{|a_{\gamma,M}|^2}+\sum_{\gamma\in C_{MN}} \log \frac{|N|}{|a_{\gamma,N}|^2}=\psi(N)S_M+\psi(M)S_N.
    \end{align*}}
    It is thus enough to first consider the case $N=P^r$ where $P$ is a monic prime polynomial and $r$ is a natural number, in which case
    {\allowdisplaybreaks
    \begin{align*}
        S_{P^r} & = \sum_{\gamma\in C_{P^r}}\log \frac{|P^r|}{|a_\gamma|^2}=\psi(P^r)\log |P^r|-2\sum_{\gamma\in C_{P^r}}\log |a_\gamma|\\
        & = \psi(P^r)\log |P^r|- 2\sum_{i=1}^r\varphi(P^{r-i})\log|P|^i\\
        & = \psi(P^r)\log |P^r|- 2\log |P|\left[\sum_{i=1}^{r-1} i|P|^{r-i}\left(1-\frac{1}{|P|}\right) +r\right]\\
        & = \left[|P|^{r-1}(|P|+1)r-\frac{2(|P^r|-1)}{|P|-1}\right]\log|P|.
    \end{align*}}
    This proves Part 1, and Part 2 follows from $|z_\gamma|_\ti = \frac{|a_\gamma|}{|d_\gamma|}|z|_\ti$.

    To prove Part 3, let $\gamma = \begin{pmatrix}
        a & b \\ 0 & d
    \end{pmatrix}$. Then we have
    \[
    D(\Lambda_{z_\gamma}) = D\left(\frac{1}{d}((az+b)A + dA)\right)
    = \frac{1}{|d|^2}D((az+b)A + dA)
    = \frac{1}{|d|^2}\det(\gamma)D(\Lambda_z)
    = \frac{|a|}{|d|}D(\Lambda_z),
    \]
    by \cite[Prop. 4.4]{Tag}. The result now follows from part 1.
\end{proof}

\begin{lemma}\label{DeltaInvariance}
    If $z\in\Omega$ and $z'\in\Omega$ are $\Gamma$-equivalent, then
    \[
    \log|\Delta(z')| = \log|\Delta(z)| - \frac{q^2-1}{2}\big(\log|z'|_{\mathrm{i}} - \log|z|_{\mathrm{i}}\big).
    \]
\end{lemma}

\begin{proof}
    Since $\Delta$ is a Drinfeld modular form of weight $q^2-1$, and from the functional equation (\ref{imtrans}) for $|\cdot|_\ti$, it follows that the map
    \[
    z \longmapsto \log\big(|\Delta(z)||z|_\ti^{\frac{q^2-1}{2}}\big)
    \]
    is $\Gamma$-invariant. The result follows.
\end{proof}

For each $\gamma\in C_N$, we write $\widetilde{z_\gamma}\in\CF$ for a fundamental representative of $z_\gamma$. Then we compute
{\allowdisplaybreaks
\begin{align}\nonumber
    h(\Phi_N(X,j(z))) =& \sum_{\gamma\in C_N}\log\max \{|j(\widetilde{z_\gamma})|,1\}\\\nonumber
    = & \sum_{\gamma\in C_N}\log\max \{|g(\widetilde{z_\gamma})|^{q+1},|\Delta(\widetilde{z_\gamma})|\}-\sum_{\gamma\in C_N}\log|\Delta(\widetilde{z_\gamma})|\\\nonumber
    = & \sum_{\gamma\in C_N}\log\max \{|g(\widetilde{z_\gamma})|^{q+1},|\Delta(\widetilde{z_\gamma})|\} - 
    \sum_{\gamma\in C_N}\log|\Delta(z_{\gamma})| \\\nonumber
    & + \frac{q^2-1}{2}\sum_{\gamma\in C_N} \big( \log|\widetilde{z_\gamma}|_\ti - \log|z_\gamma|_\ti\big),
    \;\text{by Lemma \ref{DeltaInvariance},} \\\nonumber
    = &\sum_{\gamma\in C_N}\log\max \{|g(\widetilde{z_\gamma})|^{q+1},|\Delta(\widetilde{z_\gamma})|\}- \psi(N)\log|\Delta(z)| \\ \nonumber
    & + \frac{q^2-1}{2}\sum_{\gamma\in C_N} \big( \log|\widetilde{z_\gamma}|_\ti - \log|z_\gamma|_\ti\big),\;\text{by Proposition \ref{davg},}\\\nonumber
    = &\sum_{\gamma\in C_N}\log\max \{|g(\widetilde{z_\gamma})|^{q+1},|\Delta(\widetilde{z_\gamma})|\}- \psi(N)\log|\Delta(z)| \\\nonumber 
    & +\frac{q^2-1}{2}\psi(N)\big(\deg N - 2\lambda_N - \log|z|_\ti\big)
    + \frac{q^2-1}{2}\sum_{\gamma\in C_N}\log|\widetilde{z_\gamma}|_\ti,
    \;\text{by Lemma \ref{imagsum},}\\\nonumber
    = & \frac{q^2-1}{2}\psi(N)\big(\deg N - 2\lambda_N\big)
    - \psi(N)\log\left[|\Delta(z)||z|_\ti^{\frac{q^2-1}{2}}\right]\\ \label{Mahler}
    & + \sum_{\gamma\in C_N}\log\max \{|g(\widetilde{z_\gamma})|^{q+1},|\Delta(\widetilde{z_\gamma})|\}
    +\frac{q^2-1}{2}\sum_{\gamma\in C_N}\log|\widetilde{z_\gamma}|_\ti.
\end{align}}

Since the $\widetilde{z_\gamma}$ are all reduced, the first sum above can be estimated using \cite[Lemma 5.3]{BPR}, which says that, for $z$ reduced,
\begin{equation}\label{Lemma5.3}
    q^2 \leq \log\max\{ |g(z)|^{q+1}, |\Delta(z)|\} \leq q^2 + q.
\end{equation}
In particular, we have
\begin{equation}\label{eq:analyticbound}
    q^2\psi(N) \leq \sum_{\gamma\in C_N}\log\max 
    \{|g(\widetilde{z_\gamma})|^{q+1},|\Delta(\widetilde{z_\gamma})|\}
    \leq (q^2 + q)\psi(N).
\end{equation}

Our main task now is to estimate the second sum, which we split into one sum with small $d_\gamma$ and one with large $d_\gamma$:
\[
\sum_{\gamma\in C_N}\log|\widetilde{z_\gamma}|_\ti = 
\sum_{\substack{|d_{\gamma}|>\sqrt{|N||z|_{\ti}}\\ \gamma \in C_N}}
\log|\widetilde{z_\gamma}|_{\ti}
+\sum_{\substack{|d_{\gamma}|\leq \sqrt{|N||z|_{\ti}}\\ \gamma \in C_N}}
\log|\widetilde{z_\gamma}|_\ti,
\]
where $\gamma=\begin{pmatrix}
    a_\gamma & b_{\gamma}\\
    0 & d_{\gamma}
\end{pmatrix}\in C_N$. 

\subsection{The case $|d_\gamma| \leq \sqrt{|N||z|_{\mathrm{i}}}$}

\begin{prop}\label{smalld}
Suppose that $|z|_{\mathrm{i}} \geq 1$. Then we have
    \[
    \sum_{\substack{|d_{\gamma}|\leq \sqrt{|N||z|_{\mathrm{i}}}\\ \gamma \in C_N}}
\log|\widetilde{z_\gamma}|_{\mathrm{i}} \leq 
\psi(N)\left(\frac{1}{q} + \log|z|_{\mathrm{i}}\right).
    \]
\end{prop}

\begin{proof}
    When $|d_\gamma| \leq \sqrt{|N||z|_\ti}$, we find that
\[
|z_{\gamma}|_\ti = \left|\frac{a_{\gamma}}{d_{\gamma}}\right|\cdot |z|_\ti
= \frac{|N|}{|d_{\gamma}|^2}\cdot |z|_\ti
\geq 1,
\]
so by Lemma \ref{CFimag}  $|\widetilde{z_\gamma}|_\ti = |z_\gamma|_\ti$. 
Now we compute (writing $a_\gamma, b_\gamma, d_\gamma$ in $a,b,d$ respectively for short, and $e_d=\gcd(d,N/d)$): % (where we denote $r=\gcd(a,d) = \gcd(N/a, d)$):
{\allowdisplaybreaks
\begin{align*}
    \sum_{\substack{|d|\leq \sqrt{|N||z|_{\ti}}\\ \gamma \in C_N}}
\log|\widetilde{z_\gamma}|_\ti 
& = \sum_{\substack{d|N\\ |d|\leq \sqrt{|N||z|_{\ti}}}} 
\sum_{\substack{|b|<|d|\\ \gcd(b,e_d)=1}} \log\left(\frac{|a|}{|d|}|z|_\ti\right)\\
& = \sum_{\substack{d|N\\ |d|\leq \sqrt{|N||z|_{\ti}}}} \frac{|d|\varphi(e_d)}{|e_d|} \left[\log|z|_\ti + \log\frac{|a|}{|d|}\right] \\
& \leq \psi(N)\log|z|_\ti + \sum_{\substack{d|N\\ |d|\leq \sqrt{|N||z|_{\ti}}}} \frac{|a|\varphi(e_d)}{q|e_d|}, \;\text{since}\; \log\frac{|a|}{|d|} \leq \frac{1}{q}\frac{|a|}{|d|},\\
& \leq  \psi(N)\log|z|_\ti + \frac{1}{q}\sum_{a|N} \frac{|a|\varphi(e_d)}{|e_d|}\\
& = \psi(N)\left(\frac{1}{q} + \log|z|_{\ti}\right).
\end{align*}}
\end{proof}

\subsection{The case $|d_\gamma|>\sqrt{|N||z|_\mathrm{i}}$} We take $z\in \mathcal{F}$ so that $|z|=|z|_\ti\geq 1$.

The difficulty is that, while we have an exact expression for $z_\gamma$, we do not have good control over the fundamental representative $\widetilde{z_\gamma}$. Instead, we will use the 
function field analogues of Farey sequences worked out by L.~C.~Hsia in \cite[Section 3]{hsia}. 
This allows us to approximate each $\widetilde{z_\gamma}$ by an explicit equivalent point $\widehat{z_\gamma}$, which we prove is close to being reduced. The needed material has been covered in paragraph \ref{Farey} and paragraph \ref{Funda}, so we now have all the ingredients needed to estimate the sum over large $d_\gamma$.

\begin{prop}\label{bigd}
    Suppose that $|z|_{\mathrm{i}} \geq 1$. Then we have
    \[
    \sum_{\substack{|d_{\gamma}|> \sqrt{|N||z|_\mathrm{i}}\\ \gamma \in C_N}}
\log|\widetilde{z_\gamma}|_\mathrm{i} \leq
\psi(N)\left(4 + \frac{2q}{(q-1)^2}\right).
    \]
\end{prop}

\begin{proof}
    For each $\gamma\in C_N$ with $d_\gamma > \sqrt{|N||z|_{\ti}}$, denote by $\widehat{z_\gamma}\in\Omega$ the point constructed in Section \ref{Funda}, which is $\Gamma$-equivalent to $\widetilde{z_\gamma}$. 
    By Lemma \ref{zhat} we have $\log|\widetilde{z_\gamma}|_{\ti}\leq \log|\widehat{z_\gamma}|_{\ti}+4$, so it remains to estimate the following sum. We set $y=|z|_\ti$. As before, we write $a,b,d$ for $a_\gamma, b_\gamma$ and $d_\gamma$, respectively. Set again $M:=\left\lceil \log\frac{|d|}{\sqrt{|N||z|_\ti}}\right\rceil\geq 1.$ We have:
{\allowdisplaybreaks
\begin{align*}
    S := & \sum_{\substack{|d|> \sqrt{|N|y}\\ \gamma \in C_N}}
\log|\widehat{z_\gamma}|_\ti \\
     \leq & \sum_{\substack{|d|> \sqrt{|N|y}\\ \gamma \in C_N}}
     \log\frac{|d|^2}{|N|y|f|^2} \\
     = & \sum_{\substack{|d|> \sqrt{|N|y}\\ \gamma \in C_N}} 
     \sum_{\substack{\deg f < M\\ \text{$f$ monic}}}
     \underbrace{\sum_{\substack{|h| < |f|\\ \gcd(h,f)=1}}
     \sum_{\substack{\frac{b}{d}\in D_M(h/f)\\ \gcd(b,e_d)=1}}}_{(*)}
     \big(2\log|d| - \log(|N|y) - 2\log|f|\big),
\end{align*}}

\noindent where in the four summations we have that only $d$ runs in the first one (not $\gamma$), $f$
runs in the second, $h$ runs in the third, and $b$ runs in the fourth. 
First, by Lemma \ref{countD_M}, for each $d$ and $f$, 
the total number of summands in (*) is 
\[
\frac{\varphi(f)}{|f|}\frac{\varphi(e_d)}{|e_d|}\frac{|d|}{q^M}.
\]

Next, we compute the elementary sum
{\allowdisplaybreaks
\begin{align*}
    \sum_{\substack{\deg f < M\\ \text{$f$ monic}}} \log|f| 
    & = \sum_{n=1}^{M-1} nq^n = q\cdot \frac{\mathrm{d}}{\mathrm{d}q} \sum_{n=1}^{M-1} q^n \\
    & = q\cdot \frac{\mathrm{d}}{\mathrm{d}q} \sum_{n=0}^{M-1} q^n = 
    q\cdot \frac{\mathrm{d}}{\mathrm{d}q} \left[\frac{q^{M}-1}{q-1}\right] \\
    & = \frac{(q-1)Mq^M - q^{M+1} +q}{(q-1)^2}.
\end{align*}}

We may thus estimate our sum as follows.
{\allowdisplaybreaks
\begin{align*}
    S & \leq 
    \sum_{\substack{|d|> \sqrt{|N|y}\\ \gamma \in C_N}} 
    \frac{\varphi(e_d)}{|e_d|}\frac{|d|}{q^M}
    \sum_{\substack{\deg f < M\\ \text{$f$ monic}}} \frac{\varphi(f)}{|f|} \big(2\log|d| - \log(|N|y) - 2\log|f|\big) \\
    & \leq 
    \sum_{\substack{|d|> \sqrt{|N|y}\\ \gamma \in C_N}} 
    \frac{\varphi(e_d)}{|e_d|}\frac{|d|}{q^M}
    \sum_{\substack{\deg f < M\\ \text{$f$ monic}}} \big(2\log|d| - \log(|N|y) - 2\log|f|\big) \\
    & \leq 
    \sum_{\substack{|d|> \sqrt{|N|y}\\ \gamma \in C_N}} 
    \frac{\varphi(e_d)}{|e_d|}\frac{|d|}{q^M}
    \sum_{\substack{\deg f < M\\ \text{$f$ monic}}} \big(2M - 2\log|f|\big) \\
    & =
    \sum_{\substack{|d|> \sqrt{|N|y}\\ \gamma \in C_N}} 
    \frac{\varphi(e_d)|d|}{|e_d|} \frac{2}{q^M}\left( 
    M\frac{q^M-1}{q-1} - \frac{(q-1)Mq^M - q^{M+1} +q}{(q-1)^2}
    \right)\\
    & = 
    \sum_{\substack{|d|> \sqrt{|N|y}\\ \gamma \in C_N}} 
    \frac{\varphi(e_d)|d|}{|e_d|} \frac{2}{q-1}\left( 
    \frac{q}{q-1} - \frac{M}{q^M} -\frac{q}{q^M(q-1)}
    \right)\\ 
    & \leq 
    \sum_{\substack{|d|> \sqrt{|N|y}\\ \gamma \in C_N}} 
    \frac{\varphi(e_d)|d|}{|e_d|} \frac{2q}{(q-1)^2} \\
    & \leq \psi(N) \frac{2q}{(q-1)^2}.
\end{align*}}
\end{proof}

\subsection{Specialization}

We now determine the height of $\Phi_N(X,Y)$ from the height of a suitable specialization $\Phi_N(X,j)$. In characteristic zero, this is done via Lagrange interpolation and a number of carefully chosen $j$-values. Over function fields, however, it turns out that specialization at a single suitable value of $j$ suffices.

\begin{lemma}\label{specialization}
    Let $f(X,Y)\in \Cinf[X,Y]$ be a non-zero polynomial. 
    Then there exists $y\in\overline{\BF_q}$ such that
    $h(f(X,Y)) = h(f(X,y))$.
\end{lemma}

% [CHECK: should we restrict the lemma to $F_\infty$, and maybe add a remark concerning $\mathbb{C}_\infty$ to avoid the referee's request for a full proof?]

\begin{proof}
    First, we treat the case where the coefficients of $f(X,Y)$ are in $F_\infty = \BF_q((\frac{1}{t}))$, as is the case when $f(X,Y) = \Phi_N(
    X,Y)$, and this case is simpler.
    The height $h(f(X,Y))=m$ equals the largest exponent of $t$ appearing in the coefficients of~$f$.
    We write 
    \begin{align*}
        f(X,Y) & = \sum_{i=0}^d a_i(t,Y) X^i, \qquad a_i(t,Y)\in\BF_q((\tfrac{1}{t}))[Y]\\
        & = \sum_{i=0}^d X^i \sum_{j=-\infty}^{m}b_{i,j}(Y)t^j, \qquad b_{i,j}\in\BF_q[Y].
    \end{align*}
    Now let $y\in\overline{\BF_q}$ be any point which is not a root of any of the non-zero $b_{i,m}\in\BF_q[Y]$. Then, since each non-zero $|b_{i,j}(y)|=1$, the height of $f(X,y)$ again equals $m$.

    In the general case, it is known that the elements of $\Cinf$ can be represented as generalized power series (also known as {\em Hahn-Mal'cev-Neumann series}) of the form $\displaystyle{x = \sum_{k\in K}x_kt^{-k}}$, where $K\subset\BQ$ is a well-ordered subset depending on $x$, and each $x_k\in\overline{\BF_q}$; see \cite[Cor. 11]{Ked01} and \cite[\S15]{Ked17}. 
    We adopt the convention where $x_{-m}\neq 0$ for $m=-\min K$, so $|x| = q^m$.

    We write
    \[
    f(X,Y) = \sum_{i=0}^d X^i \sum_{j=0}^d \left( \sum_{k\in K_{ij}} x_{i,j,k}t^{-k}\right)Y^j,
    \]
    where each $K_{ij}\subset\BQ$ is a well-ordered set and each $x_{i,j,k}\in\overline{\BF_q}$. 
    Let $m = -\min \bigcup_{i,j} K_{i,j}$, which exists by well-ordering. 
    By our convention, $x_{i,j,-m}\neq 0$ for at least one pair $(i,j)$. It follows that 
    $m = h(f(X,Y))$ is the largest power of $t$ appearing in $f(X,Y)$. We also set $n = -\min \left(\bigcup_{i,j} K_{i,j} \smallsetminus \{-m\}\right) < m$, which again exists by well-ordering.
    Now we decompose
    \[
    f(X,Y) = \sum_{i=0}^d X^i \sum_{j=0}^d \left( \sum_{\substack{k\in K_{ij} \\ k \geq -n}} x_{i,j,k}t^{-k}\right)Y^j
    \; + \;  \sum_{i=0}^d X^i \underbrace{\sum_{j=0}^d x_{i,j,-m}t^m Y^j}_{b_i(Y)t^m}.
    \]
    As before, at least one $b_i(Y) = \sum_{j=0}^d x_{i,j,-m}Y^j \in \overline{\BF_q}[Y]$ is non-zero, and we pick $y\in\overline{\BF_q}$ such that $b_j(y)\neq 0$, and thus $|b_j(y)|=1$, for every non-zero $b_j(Y)$.
    Once again, the largest power of $t$ appearing in $f(X,y)$ is $m=h(X,Y)$.
    % In the general case, it is known that the elements of $\Cinf$ can be represented as generalized power series (also known as {\em Hahn-Mal'cev-Neumann series}) of the form $\displaystyle{x = \sum_{j\in J}x_jt^{-j}}$, where $J\subset\BQ$ is a well-ordered subset depending on $x$, see \cite[Cor. 11]{Ked01} and \cite[\S15]{Ked17}. The above argument allows us to conclude in this case, too. [CHECK: DETAILS!!! :) referee is asking]
\end{proof}

We easily obtain the following generalization, although we won't need it.

\begin{lemma}
    Let $n\geq 2$ and $f(X_1,X_2,\ldots,X_n) \in \Cinf[X_1,X_2,\ldots,X_n]$ be non-zero. 
    Let $1 < m \leq n$. Then there exist $x_m, x_{m+1}, \ldots, x_n \in \overline{\BF_q}$ such that
    $h(f(X_1,X_2,\ldots,X_n)) = h(f(X_1, X_2, \ldots, X_{m-1}, x_m, x_{m+1} \ldots, x_n))$.
\end{lemma}

\begin{proof}
    By induction, it suffices to prove this for $m=n$. Now we follow the general part of the proof of Lemma \ref{specialization}, except that the index $i$ is replaced by the multi-index $(i_1,i_2,\ldots,i_{n-1})$. 
    % For instance,
    % each appearance of $\sum_{i=0}^dX^i$ is replaced by $\sum_{i_1, i_2, \ldots, i_{n-1}} X_1^{i_1}X_2^{i_2}\cdots X_{n-1}^{i_{n-1}}$.
\end{proof}

We are now ready to prove our main result. By Lemma \ref{specialization}, there exists  $y=j(z)\in\overline{\BF_q}$ for which
\[
h(\Phi_N(X,Y)) = h(\Phi_N(X,j(z))) = \sum_{\gamma\in C_N}\log\max(1, |j(z_\gamma)|).
\]
We may choose $z\in\CF$. Since $y=j(z)\in\overline{\BF_q}$ and thus $|j(z)|=1$, this imposes extra conditions on $z$.

\begin{lemma}
    Suppose that $z\in\CF$ and $\log|j(z)|\leq q$. Then $|z|_{\mathrm{i}} = 1$ and $\log|\Delta(z)|=q^2$.
\end{lemma}

\begin{proof}
    Given that $z\in\CF$ and $\log|j(z)|\leq q$, \cite[Lemma 5.4]{BPR} implies that $|z|_\ti=1$.
      Now \cite[Theorem 2.13]{Gek0} gives us $\log|\Delta(z)| = q^2$.
\end{proof}

Putting everything together, with $|z|_\ti=1$ and $\log|\Delta(z)|=q^2$, we compute:
{\allowdisplaybreaks
\begin{align*}
    h(\Phi_N(X,Y)) = &  h(\Phi_N(X,j(z)) = \sum_{\gamma\in C_N}\log\max(|j(\widetilde{z_\gamma})|, 1) \\
    = & \frac{q^2-1}{2}\psi(N)\big(\deg N - 2\lambda_N\big)
    - \psi(N)\log\left[|\Delta(z)||z|_\ti^{\frac{q^2-1}{2}}\right]\\
    & + \sum_{\gamma\in C_N}\log\max \{|g(\widetilde{z_\gamma})|^{q+1},|\Delta(\widetilde{z_\gamma})|\}
    +\frac{q^2-1}{2}\sum_{\gamma\in C_N}\log|\widetilde{z_\gamma}|_\ti \\
    = & \frac{q^2-1}{2}\psi(N)\left(\deg N - 2\lambda_N - \frac{2q^2}{q^2-1} \right) \\
    & + \sum_{\gamma\in C_N}\log\max \{|g(\widetilde{z_\gamma})|^{q+1},|\Delta(\widetilde{z_\gamma})|\}
    +\frac{q^2-1}{2}\sum_{\gamma\in C_N}\log|\widetilde{z_\gamma}|_\ti.
\end{align*}}

Now the upper bounds from (\ref{eq:analyticbound}), Proposition \ref{smalld} and Proposition \ref{bigd} give us
{\allowdisplaybreaks
\begin{align*}
    h(\Phi_N(X,Y)) \leq &  \frac{q^2-1}{2}\psi(N)\left(\deg N - 2\lambda_N - \frac{2q^2}{q^2-1} \right) \\
    & + \psi(N)(q^2+q) + \frac{q^2-1}{2}\psi(N)\left(4 + \frac{2q}{(q-1)^2} + \frac{1}{q}\right) \\
    = & \frac{q^2-1}{2}\psi(N)\left(\deg N - 2\lambda_N + 4 + \frac{5q^3-q^2-q+1}{q(q+1)(q-1)^2}
    \right).
\end{align*}}

On the other hand, combining the lower bound from (\ref{eq:analyticbound}) with the trivial $|\widetilde{z_\gamma}|_\ti \geq 1$, we obtain
\begin{align*}
 h(\Phi_N(X,Y)) \geq &  \frac{q^2-1}{2}\psi(N)\left(\deg N - 2\lambda_N - \frac{2q^2}{q^2-1} \right) 
     + \psi(N)q^2  \\
    = & \frac{q^2-1}{2}\psi(N)(\deg N - 2\lambda_N).
\end{align*}
This completes the proof of Theorem \ref{thm:main}. \qed

\section{Heights of Hecke images}\label{proof of Hecke}

Let $K$ be a finite extension of $F=\BF_q(t)$. Throughout this section, we assume $[K:F]$ is large enough so that all Drinfeld modules considered are defined over $K$, and are of rank 2.

Let $\phi$ be a Drinfeld $A$-module of rank 2 defined over $K$ and for every submodule $C\subset \phi[N]$ with $C\cong A/NA$ we denote by $\phi/C$ the quotient Drinfeld module, \textit{i.e.} the image under the isogeny $f : \phi \rightarrow \phi/C$ with $\ker f = C$. 

In this section, we give a proof of Theorem \ref{Silverman}, which gives a comparison of the height of $\phi$ to the average of the heights of $\phi/C$ for all $C\cong A/NA$ (the case of a single isogeny was studied in \cite{BPR}).

\subsection{Proof of Theorem \ref{Silverman}}

The strategy to prove Theorem \ref{Silverman} is to split $h(j(\phi))$ and $h(j(\phi/C))$ into a finite part and an infinite part, and then compute the difference separately.

Consider a Drinfeld module $\rho$ of rank 1 over an $A$-field $L$, and $N\in A$ be monic away from the characteristic. Denote the action of $A$ on $L$ through $\rho$ by $*$. Take any $0\neq \omega\in L$ and view $\frac{L}{A*\omega}$ as an $A$-module through $\rho$. 

\begin{lemma}\label{Nsubmodule}
    If $\rho_N(X)-\omega$ (resp. $\rho_N(X)$) splits in $L$ such that for some $x,y\in L$, $\rho_N(x)=\omega$ (resp. $\rho_N(y)=0$) with $\rho_a(x)\neq \omega$ (resp. $\rho_a(y)\neq 0$) for all $a|N$, $\deg(a) < \deg(N)$, then the $A$-submodules of $\frac{L}{A*\omega}$ of order $|A/N|$ and annihilated by $N$ are given by
    \[
    \frac{A(a*x+b*y)+(d)*y}{A*\omega},
    \]
    where $a,b,d\in A$, $a,d$ monic with $ad=N$, and $\deg(b)<\deg(d)$.
\end{lemma}

\begin{proof}
    It is easy to verify that each of these submodules are annihilated by $N$ and are of order $|A/N|$. Let $g\in L$ such that $\rho_N(g)\in A*\omega$. It is therefore of the form $g=a*x+b*y$ for some $a,b\in A$. Taking the quotient gives us $\deg(a),\deg(b)<\deg(N)$.
    
    Let $a_1:=\gcd(a,N)$, $a_2:=a/\gcd(a,N)$, $d_1=N/a_1$ and $b_1\in A$ such that $b_1\equiv a_2^{-1}b\ (\text{mod}\ d_1)$. Hence we have
    \[
    g\in \frac{A(a_1*x+b_1*y)+(d_1)*y}{A*\omega},
    \]
    One can then easily deduce our desired result.
\end{proof}

\begin{prop}\label{finite places}
    Let $\phi$ be a rank 2 Drinfeld module over $\bar{F}$ and $v$ be a finite place. Then
    \[
    \log\max(1, |j(\phi)|_v) = \frac{1}{\sigma_1(N)}\sum_C \log\max(1, |j(\phi/C)|_v),
    \]
    where $C$ runs through all $A$-submodules of $\phi[N]$ of order $|A/N|$ and annihilated by $N$, and $\sigma_1(N)=\sum\limits_{\substack{a|N,\\ a\ monic}}|a|.$
\end{prop}

\begin{proof}
If $v(j(\phi))\geq 0$, then $\phi$ has good reduction at $v$. Since $\phi$ is isogenous to $\phi/C$, the module $\phi/C$ has good reduction at $v$ as well. In this case, both sides vanish. 

Assume $v(j(\phi))<0$. Then $\phi$ has stable bad reduction at $v$. Let $\bar{F}_v$ be the completion of $\bar{F}$ at $v$, and $\Ocal$ be the valuation ring of $\bar{F}_v$. Applying the Tate uniformisation, we obtain a pair $(\rho,\Lambda)$ corresponding to $\phi$ where $\rho$ is a Drinfeld module of rank $1$ over $\Ocal$ and $\Lambda\subset \bar{F}_v$ is an $A$-lattice of rank 1 with the action given by $\rho$. Write $\Lambda=A*\omega$, where $*$ is the action of $A$ through $\rho$ on $\bar{F}_v$. Taking isomorphism if necessary, we may assume $\rho_t=t+\tau$. Moreover, we have an isomorphism of $A$-modules
\[
\frac{\bar{F}_v}{A*\omega}\simeq \bar{F}_v,
\]
where the left-hand side is an $A$-module via $\rho$ while the right-hand side has the $A$-module structure by $\phi$. By the calculation of \cite[Example 6.2.4]{Papikian} we see
\begin{equation}\label{valjinv}
    v(j(\phi))=(q-1)v(\omega)<0.
\end{equation}
 Take $x\in \bar{F}_v$ (resp. $y\in \bar{F}_v$) such that $\rho_N(x)=\omega$ (resp. $\rho_N(y)=0$) with $\rho_a(x)\neq \omega$ (resp. $\rho_a(y)\neq 0$) for all $a|N$, $\deg(a)<\deg(N)$. Lemma \ref{Nsubmodule} tells us that all the $A$-submodules of $\frac{\bar{F}_v}{A*\omega}$ of order $|A/N|$ and annihilated by $N$ are given by
\[
    \frac{A(a*x+b*y)+(d)*y}{A*\omega},
    \]
where $a,b,d\in A$, $a,d$ monic with $ad=N$, and $\deg(b)<\deg(d)$. Since $\rho_N(x)=\omega$ and $\rho_N(y)=0$, one can easily check that $v(x)<0$ and $v(y)\geq 0$. We write the quotient Drinfeld modules using the isomorphism as follows:
\[
\frac{\bar{F}_v}{A(a*x+b*y)+(d)*y}\xrightarrow{\sim} \frac{\bar{F}_v}{A(a^2*x+ab*y)},\quad \xi\mapsto a*\xi. 
\]
We observe that $v(ab*y)\geq 0$ and $v(a^2*x)=|a|^2v(x)<0$, which implies $v(a^2*x+ab*y)=|a|^2v(x)$. Again, oberve that $v(\omega)=|N|v(x)$ so that we get 
\[
v(a^2*x+ab*y)=\frac{|a|}{|d|}v(\omega).
\]
Applying (\ref{valjinv}) to the quotient Drinfeld modules we see $v(j(\phi/C))$ is given by some $v(a^2*x+ab*y)$. Summing over all submodules $C$ gives
\[
\sum_{C} v(j(\phi/C)) = \sum_C (q-1)\frac{|a|}{|d|}v(\omega) = (q-1)v(\omega)\sum_{d|N}|d|\frac{|a|}{|d|} = \sigma_1(N)v(j(\phi)).
\]
\end{proof}

For our purposes, we only want to sum over those submodules $C$ which are cyclic. This leads to the additional condition on the triple $(a,b,d)$ that $\gcd(a,b,d)=1$. With this restriction, in the last equation of the above proof becomes
\[
\sum_{\text{$C$ cyclic}} v(j(\phi/C)) = \sum_{\scriptsize \begin{pmatrix}
    a & \!\!\!\!\!b \\ 0 & \!\!\!\!\!d
\end{pmatrix}\in C_N} (q-1)\frac{|a|}{|d|}v(\omega) = (q-1)v(\omega)\sum_{d|N}\frac{\varphi(e_d)|d|}{|e_d|}\frac{|a|}{|d|} = \psi(N)v(j(\phi)),
\]
where $e_d = \gcd(a,d) = \gcd(N/d,d)$. 
We thus obtain the following corollary:

\begin{corollary}\label{cyclic submodules}
        Let $\phi$ be a rank 2 Drinfeld module over $\bar{F}$ and $v$ be a finite place. Then
    $$\log\max(1, |j(\phi)|_v) = \frac{1}{\psi(N)}\sum_{C} \log\max(1, |j(\phi/C)|_v),$$
    where $C$ runs through all cyclic $A$-submodules of $\phi[N]$ of order $|A/N|$.
\end{corollary}

\begin{proof}[Proof of Theorem \ref{Silverman}]
    By Corollary \ref{cyclic submodules}, we have
    \begin{align*}
        & h(j(\phi)) - \frac{1}{\psi(N)}\sum_{C}h(j(\phi/C))  = 
        h^\infty(j(\phi)) - \frac{1}{\psi(N)}\sum_{C}h^\infty(j(\phi/C)) \\
        & = \frac{1}{[K:F]}\sum_{\sigma : K \hookrightarrow\Cinf} n_\sigma
        \left[\log\max\big(1, |\sigma(j(\phi))|\big) - \frac{1}{\psi(N)} 
        \underbrace{\sum_C\log\max\big(1, |\sigma(j(\phi/C))|\big)}_{(*)} 
        \right].
    \end{align*}
    Now by Lemma \ref{mahler}  the sum $(*)$ equals the height of the specialized modular polynomial $\Phi_N\big(X,\sigma(j(\phi))\big)$, which is easily seen to be bounded as follows
    \[
        h\big(\Phi_N(X, \sigma(j(\phi)))\big) \leq h(\Phi_N) + \psi(N)\log\max\big(1, |\sigma(j(\phi))|\big).
    \]
    We thus obtain a lower bound of $-h(\Phi_N)/\psi(N)$, and the lower bound in
    Theorem \ref{Silverman} now follows from Theorem \ref{thm:main}.

    Next, we prove the upper bound. Fix $\sigma : K \hookrightarrow \Cinf$ and write $\sigma(j(\phi)) = j(z)$ for some $z = z_\sigma\in\CF$. We rewrite the term in the square brackets using (\ref{Mahler}) as follows.
    {\allowdisplaybreaks
    \begin{align*}
        & \log\max(1, |j(z)|) - \frac{1}{\psi(N)}\sum_{\gamma\in C_N}\log\max(1, |j(z_\gamma)|) \\
        = & \;\big(\log\max(|g(z)|^{q+1}, |\Delta(z)|) - \log|\Delta(z)|\big) 
        -\frac{q^2-1}{2}(\deg N - 2\lambda_N) + \log\left[|\Delta(z)||z|_\ti^{\frac{q^2-1}{2}}\right] \\
        & - \frac{1}{\psi(N)} \sum_{\gamma\in C_N}\log\max 
        \{|g(\widetilde{z_\gamma})|^{q+1},|\Delta(\widetilde{z_\gamma})|\}
        - \frac{q^2-1}{2}\frac{1}{\psi(N)}\sum_{\gamma\in C_N}\log|\widetilde{z_\gamma}|_\ti \\
        = & \underbrace{\left[ \log\max(|g(z)|^{q+1}, |\Delta(z)|)  
        - \frac{1}{\psi(N)} \sum_{\gamma\in C_N}\log\max 
            \{|g(\widetilde{z_\gamma})|^{q+1},|\Delta(\widetilde{z_\gamma})|\}\right]}_{(**)} \\
        & -\frac{q^2-1}{2}(\deg N - 2\lambda_N) + 
        \frac{q^2-1}{2}\underbrace{\left[\log|z|_\ti - \frac{1}{\psi(N)}\sum_{\gamma\in C_N}\log|\widetilde{z_\gamma}|_\ti
        \right]}_{(***)}.
    \end{align*}}
    By (\ref{Lemma5.3}) and (\ref{eq:analyticbound}) we have $(**) \leq q$, so it remains to bound $(***)$ from above.

    Since $|\widetilde{z_\gamma}|_\ti \geq \max(|z_\gamma|_\ti, 1)$, by Lemma \ref{CFimag}, we find that 
    \[
    \frac{1}{\psi(N)}\sum_{\gamma\in C_N}\log|\widetilde{z_\gamma}|_\ti \geq
    \max\left( 0, \frac{1}{\psi(N)}\sum_{\gamma\in C_N} \log|z_\gamma|_\ti
    \right) =
    \max(0, \log|z|_\ti-\deg N + 2\lambda_N),
    \]
    by Lemma \ref{imagsum}. Next, \cite[Lemma 5.4]{BPR} gives us 
    \[
    \log|z|_\ti \leq \log \max\left(\frac{1}{q}\log|j(z)|, 1\right).
    \]
    Now we put everything together:
    {\allowdisplaybreaks
    \begin{align*}
         & h(j(\phi)) - \frac{1}{\psi(N)}\sum_{C}h(j(\phi/C))  \\
         & =  \frac{1}{[K:F]}\sum_{\sigma : K \hookrightarrow\Cinf} n_\sigma
        \left[\log\max\big(1, |\sigma(j(\phi))|\big) - \frac{1}{\psi(N)} 
        \sum_C\log\max\big(1, |\sigma(j(\phi/C))|\big)
        \right]\\
        & \leq \frac{1}{[K:F]}\sum_{\sigma : K \hookrightarrow\Cinf} n_\sigma
        \left[
        q - \frac{q^2-1}{2}\left( \deg N - 2\lambda_N
        - \log|z_\sigma|_\ti + \max(0, \log|z_\sigma|_i -\deg N + 2\lambda_N)\big)\right)
        \right]\\
         & = \frac{1}{[K:F]}\sum_{\sigma : K \hookrightarrow\Cinf} n_\sigma
        \left[
        q + \frac{q^2-1}{2}\min\left( \log|z_\sigma|_\ti -\deg N + 2\lambda_N, \; 0
        \right)
        \right]\\
        & \leq q +\frac{q^2-1}{2}\min \left[0, \;
        2\lambda_N - \deg N + \frac{1}{[K:F]}\sum_{\sigma: K \hookrightarrow\Cinf}n_\sigma 
        \log \max\left(\frac{1}{q}\log|j(\sigma(\phi))|, 1\right)
        \right] \\
        & = q +\frac{q^2-1}{2}\min \left[0, \;
        2\lambda_N - \deg N + \log\prod_{\sigma: K \hookrightarrow\Cinf}
        \max\left(\frac{1}{q}\log|j(\sigma(\phi))|, 1\right)^{\frac{n_\sigma}{[K:F]} }
        \right] \\
        & \leq q +\frac{q^2-1}{2}\min \left[0, \;
        2\lambda_N - \deg N + \log\frac{1}{[K:F]}\sum_{\sigma: K \hookrightarrow\Cinf} n_\sigma
        \max\left(\frac{1}{q}\log|j(\sigma(\phi))|, 1\right)
        \right] (\diamond) \\
        & \leq q +\frac{q^2-1}{2}\min \left[0, \;
        2\lambda_N - \deg N + \log\frac{1}{[K:F]}\sum_{\sigma: K \hookrightarrow\Cinf} n_\sigma
        \left(
        \log\max\left(|j(\sigma(\phi))|^{1/q}, 1\right)+1\right)
        \right] \\
         & \leq q +\frac{q^2-1}{2}\min \left[0, \;
        2\lambda_N - \deg N + \log\left(\frac{1}{q}h(j(\phi)) + 1\right)
        \right],
    \end{align*}}
    where we have used the AM-GM inequality to reach step $(\diamond)$. This completes the proof of Theorem \ref{Silverman}.
\end{proof}

\subsection{Proof of Theorem \ref{Autissier}}

Let $\phi$ be a Drinfeld module defined over $K$, and suppose that $\phi$ and $\phi/C$ have everywhere stable reduction over $K$, for all submodules $C\subset \phi[N]$ with $C\cong A/NA$. Recall that we denote by $h_\Tag(\phi)$ the stable Taguchi height of $\phi$ (see paragraph \ref{heights}, or \cite{Tag}, or \cite{Wei} for definitions). 

\begin{proof}[Proof of Theorem \ref{Autissier}]
    For each place $v\in M_K^\infty$, by the proof of \cite[Lemma 5.1]{BPR}, the quantity     $|\Delta(\phi)|_v D(\Lambda_v)^{(q^2-1)/2}$ is an isomorphism invariant. 
    Therefore, 
    we may assume that $\Lambda_v = \Lambda_{z_v}$ for some $z_v\in\Omega$ and that the Hecke images $\phi/C$ correspond to the lattices $\Lambda_{C,v} = \Lambda_{\gamma(z_v)}$ for $\gamma\in C_N$. 

    Now we compute from (\ref{eq: Tag height 2})
    {\allowdisplaybreaks
    \begin{align*}
        & h_\Tag(\phi) - \frac{1}{\psi(N)}\sum_C h_\Tag(\phi/C) \\
        = & \; \frac{1}{(q^2-1)[K:F]}\sum_{v\in M_K^f}n_v\underbrace{\left[
         \log\max(|j(\phi)|_v, 1) 
        - \frac{1}{\psi(N)}\sum_C \log\max(|j(\phi/C)|_v, 1)
        \right]}_{(*)} \\
        & - \frac{1}{(q^2-1)[K:F]}\sum_{v\in M_K^\infty}n_v\underbrace{\left[
        \log|\Delta(\phi)|_v  
        - \frac{1}{\psi(N)}\sum_C \log|\Delta(\phi/C)|_v  
        \right]}_{(**)} \\
        & - \frac{1}{2[K:F]}\sum_{v\in M_K^\infty}n_v\underbrace{\left[
        \log D(\Lambda_{v}) 
        - \frac{1}{\psi(N)}\sum_{C}\log D(\Lambda_{C, v})
        \right]}_{(***)}.
    \end{align*}}
    We have $(*)=0$ by Corollary \ref{cyclic submodules} and $(**)=0$ by Lemma \ref{davg},
    whereas $(***) = \deg N - 2\lambda_N$,  by Lemma \ref{imagsum}.3. 
    This completes the proof.

\end{proof}

\end{document}